\documentclass{gtpart}
\usepackage{amscd}
\usepackage{epsfig}
\usepackage{graphicx}  %%% the recommended graphics package
\usepackage[all]{xy}
\usepackage{amsfonts}
\title[An invariant of $3$--manifolds in $6$--manifolds and the triple linking number]{An
invariant of embeddings of $3$--manifolds in $6$--manifolds and Milnor's triple linking number}

\author{Tetsuhiro Moriyama}
\givenname{Tetsuhiro}
\surname{Moriyama}
\address{D\'epartement de Math\'ematiques\\
Institut Fourier, UMR 5582, CNRS, BP 74\\
F-38402 Saint--Martin d'H\'eres Cedex\\
France}
\email{tetsuhir@ms.u-tokyo.ac.jp}
\email{moriyama@fourier.ujf-grenoble.fr}

 %Embeddings
 %Invariants of knots and 3-manifolds
 %Isotopy
 %Knots and links in $S^3$
\arxivreference{}
\arxivpassword{}

\volumenumber{}
\issuenumber{}
\publicationyear{}
\papernumber{}
\startpage{}
\endpage{}
\doi{}
\MR{}
\Zbl{}
\received{}
\revised{}
\accepted{}
\published{}
\publishedonline{}
\proposed{}
\seconded{}
\corresponding{}
\editor{}
\version{}

\theoremstyle{plain}
\newtheorem{theorem}{Theorem}[section]    % Standard theorem environment
\newtheorem{lemma}{Lemma}[section]          % Lemma environment with numbering
\newtheorem{proposition}{Proposition}[section]          % Proposition environment with numbering
\newtheorem{corollary}{Corollary}[section]          % Corollary environment with numbering

\theoremstyle{definition}
\newtheorem{definition}{Definition}[section]    % Definition environment with
\newtheorem{remark}{Remark}[section]             % Unnumbered environment for remarks.

\makeatletter
\let\c@lemma=\c@theorem
\let\c@proposition=\c@theorem
\let\c@corollary=\c@theorem
\let\c@definition=\c@theorem
\let\c@remark=\c@theorem
\let\c@example=\c@theorem
\let\c@notation=\c@theorem
\let\c@assertion=\c@theorem
\makeatother

\theoremstyle{plain} %text of this environment is typesetted in italics
\newtheorem{mainthm}{Theorem}
\newtheorem{maincor}{Corollary}
\newtheorem{axiom}{Axiom}

\makeatletter
\let\c@maincor=\c@mainthm
\makeatother

\makeautorefname{mainthm}{Theorem}
\makeautorefname{maincor}{Corollary}
\newcommand{\N}{{\mathbb{N}}}

\newcommand{\CP}{{\mathbb{C}P}}

\newcommand{\RP}{{\RP}}
\newcommand{\Sign}{{\operatorname{Sign}\,}}

\newcommand{\bd}{\partial}
\newcommand{\HV}{{\Hat{V}}}

\newcommand{\HX}{{\Hat{X}}}

\newcommand{\HM}{{\Hat{M}}}

\newcommand{\calE}{{\mathcal{E}}}

\newcommand{\calL}{{\mathcal{L}}}

\newcommand{\calN}{{\mathcal{N}}}

\newcommand{\pair}[1]{{#1, \partial{#1}}}

\renewcommand{\a}{\alpha}
\renewcommand{\b}{\beta}
\newcommand{\g}{\gamma}
\renewcommand{\d}{\delta}

\newcommand{\s}{\sigma}
\renewcommand{\i}{\iota}

\newcommand{\bk}[1]{{\langle{#1}\rangle}}
\newcommand{\Int}{{\operatorname{Int}\,}}
\newcommand{\x}{{\times}}

\renewcommand{\L}{\Lambda}

\newcommand{\lk}{{\operatorname{lk}}}
\newcommand{\Image}{{\operatorname{Im}\,}}
\newcommand{\Ker}{{\operatorname{Ker}\,}}

\renewcommand{\O}{\Omega}
\renewcommand{\l}{\lambda}

\renewcommand{\t}{\tau}

\renewcommand{\S}{\Sigma}
\newcommand{\TX}{{\Tilde{X}}}
\renewcommand{\TH}{{\Tilde{H}}}

\newcommand{\Te}{{\Tilde{e}}}

\newcommand{\Tg}{{\Tilde{\g}}}

\newcommand{\es}{\emptyset}
\newcommand{\up}{\upsilon}

\newcommand{\Ts}{\Tilde{s}}

\newcommand{\Tc}{\Tilde{c}}

\renewcommand{\o}{\omega}

\newcommand{\Emb}{{\operatorname{Emb}}}

\newcommand{\defeq}{\stackrel{def}{\,{=}\,}}

\hypersetup{
pdftitle={An invariant of embeddings of 3-manifolds in 6-manifolds and Milnor's triple linking number},
pdfauthor={Tetsuhiro Moriyama}
}

\begin{document}
\begin{abstract}    % type your abstract below
  We give a simple axiomatic definition of a rational--valued invariant $\s(W,V,e)$ of triples
  $(W,V,e)$, where $W \supset V$ are smooth oriented closed manifolds of dimensions $6$ and $3$,
  and $e$ is a second rational cohomology class of the complement $W \setminus V$ satisfying a
  certain condition.
  The definition is stated in terms of cobordisms of such triples and the signature of
  $4$-manifolds.
  When $W = S^{6}$ and $V$ is a smoothly embedded $3$--sphere, and when 
  %$e/2$ represents the linking number with $V$, 
  $e/2$ is the Poincar\'e dual of a Seifert surface of $V$,
  the invariant coincides with $-8$ times Haefliger's embedding
  invariant of $(S^{6},V)$.
  Our definition recovers a more general invariant due to Takase, and contains a new
  definition for Milnor's triple linking number of algebraically split $3$--component links in
  $\mathbb{R}^{3}$ that is close to the one given by the perturbative series expansion of the
  Chern--Simons theory of links in $\mathbb{R}^{3}$.
\end{abstract}
\begin{asciiabstract}    % type your abstract below
  We give a simple axiomatic definition of a rational-valued invariant s(W,V,e) of triples
  (W,V,e), where W is a (smooth, oriented, closed) 6-manifold and V is a 3-submanifold of W, and
  where e is a second rational cohomology class of the complement of V satisfying a certain
  condition.
  The definition is stated in terms of cobordisms of such triples and the signature of
  4-manifolds.
  When W = S^6 and V is a smoothly embedded 3-sphere, and when e/2 is the Poincare dual of a
  Seifert surface of V, the invariant coincides with -8 times Haefliger's embedding invariant of
  (S^6,V).
  Our definition recovers a more general invariant due to Takase, and contains a new definition
  for Milnor's triple linking number of algebraically split 3-component links in R^3 that is
  close to the one given by the perturbative series expansion of the Chern-Simons theory of
  links in R^3.
\end{asciiabstract}

\maketitle

%%%%%%%%%%%%%%%%%%%%   Start of main body of article

%\makeatother

%\tableofcontents

%%%%%%%%%%%%%%%%%%%%%%%%%%%%%%%%%%%%%%%%%%%%%%%%%%%%%%%%%%%%%%%%%%%%%%%%%%%%%%%%%%%%%%%%%%%%%%%
\section{Introduction and main results}
\label{sec:intro}
%%%%%%%%%%%%%%%%%%%%%%%%%%%%%%%%%%%%%%%%%%%%%%%%%%%%%%%%%%%%%%%%%%%%%%%%%%%%%%%%%%%%%%%%%%%%%%%
Milnor~\cite{Milnor-link} proved that the link homotopy classes of algebraically split
$3$--component links $L = K_{1} \cup K_{2} \cup K_{3}$ in the Euclidean $3$--space $\R^{3}$ are
classified by the triple 
linking number $\mu(L) \in \Z$.
There are several definitions of $\mu(L)$, and one is given by
the perturbative series expression of the Chern--Simons theory of links in $\R^{3}$
(Altschuler--Freidel \cite{alt-frei}, Bar-Natan--Vassiliev \cite{barnatan-vassiliev},
Lescop \cite{lescop-links}, Thurston \cite{thurston-express}, etc.), more precisely, $\mu(L)$ is
expressed as an integral over a manifold $(T^{3} \x \R^{3}) \setminus \calL $,
where $T^{3} = S^{1} \x S^{1} \x S^{1}$ is the $3$--torus and
$\calL = T^{3}_{1} \cup T^{3}_{2} \cup T^{3}_{3} \subset T^{3} \x \R^{3}$ is the (disjoint) union of embedded $3$--tori
\[
	T^{3}_{i} = \left\{ (t_{1},t_{2}, t_{3}, x) \in T^{3} \x \R^{3}\ | \  x = f_{i}(t_{i})\right\},
\]
and where $f_{i} \co S^{1} \to \R^{3}$ is a smooth embedding representing the knot $K_{i}$.
On the other hand,
Haefliger~\cite{Haefliger-knot}~\cite{Haefliger-diff} proved that the abelian group
$\Emb(S^{3},S^{6})$ of the smooth isotopy classes of embeddings $S^3 \to S^{6}$, with the group
structure given by the connected sum, is isomorphic to $\Z$.
In this paper, we give an
interpretation of $\mu(L)$ as an invariant of the embedding $\calL \hookrightarrow T^{3} \x \R^{3}$ by
generalizing Haefliger's construction.
To this end, we will need to modify the manifold pair
$(T^{3} \x \R^{3}, \calL)$ to make it fit into our settings.
It will be replaced by $(T^{3}\x
S^{3}, M_{L})$, $M_{L} = \calL \cup (-\calL_{0})$,
so that the ambient manifold is closed and the submanifold is null--homologous, where
$S^{3} = \R^{3} \cup \left\{ \infty \right\}$,   and
$\calL_{0}$ is a $3$--submanifold of $T^{3} \x \R^{3}$ constructed from a $3$--component unlink
in $\R^{3}$ split from $L$ in the same way as we construct $\calL$ (so $M_{L}$ is the union of
$6$ disjoint copies of $T^{3}$).

In this paper, we deal with triples $\b = (Z,X,e)$, which we will call
$e$-manifolds\footnote{``$e$'' is the first letter of ``Euler class''.}, consisting of (smooth,
oriented, and compact) manifolds $Z \supset X$ of 
codimension $3$ such that $X$ is properly embedded in 
$Z$ ($\bd{X}
\subset \bd{Z}$, and $X$ is transverse to $\bd{Z}$), and a cohomology class $e \in H^{2}(Z
\setminus X;\Q)$ such that
\[
e|_{S(\nu_{X})} = e(F_{X})
\]
over $\Q$, where $e(F_{X}) \in H^{2}(S(\nu_{X});\Z)$ is the Euler class of the vertical tangent
subbundle $F_{X} \subset TS(\nu_{X})$ of the total space of the normal sphere bundle $\rho_{X} \co
S(\nu_{X}) \to X$ of $X$, and where the normal bundle $\nu_{X}$ of $X$ is identified with a tubular
neighborhood of $X$ so that we can regard $S(\nu_{X})$ as a submanifold of $Z \setminus X$.
Such a cohomology class $e$ will be called an $e$-class of $(Z,X)$ in this paper.  
%\begin{remark}
%  \label{rem:eclass}
The existence of an $e$-class implies the vanishing of the rational fundamental homology class of
$(\pair{X})$ in $(\pair{Z})$ (\fullref{prop:eclass}\,\eqref{item:quasi}), but the converse is not true in
general (\fullref{rem:S7}).
%\end{remark}

The cohomology class $e/2$ corresponds to the homology class of a Seifert surface of $X$ 
(if it exists and its normal bundle is trivial over $X$)
by the Poincar\'e duality (\fullref{cor:seifert-e}).
In particular, if $(Z,X)$ admits just one $e$-class, then $e/2$
represents the homomorphism $H_{2}(Z \setminus X;\Z) \to \Q$, $y \mapsto \lk(y, X)$, where  $\lk(y, X)$ is the
linking number of $y$ with $X$. 

Precise definitions of
$e$-class,  $e$-manifold, isomorphism (denoted by $\cong$) and cobordism of  $e$-manifolds,
etc.~will be given in \fullref{sec:preliminaries}.
We will also introduce notions of quasi $e$-class and
quasi $e$-manifold.
These are slight generalizations of $e$-class and $e$-manifold, and easier to
handle as we will explain in \fullref{rem:easier}.
We remark that, for (quasi) $e$-manifolds $\b =
(Z,X,e)$ and $\b'$ (of the same dimension), the disjoint sum $\b \amalg \b'$, the boundary $\bd{\b} =
(\bd{Z},\bd{X},e|_{\bd{Z} \setminus \bd{X}})$, and the reversing the orientation $-\b = (-Z,-X,e)$
are defined in natural ways.
By definition, an $e$-class is a quasi $e$-class, and the boundary of
a quasi $e$-manifold is an $e$-manifold (not just a quasi $e$-manifold).

There are two main purposes of the present paper. 
One is to show the existence and the uniqueness
of a rational--valued invariant $\s(\a)$ of the isomorphism classes of closed $6$--dimensional
$e$-manifolds $\a =
(W,V,e)$ ($\bd{W} = \bd{V} = \es$, $\dim W = 6$) which is uniquely characterized by the
following two axioms (\fullref{thm:main1}):
%\vspace{2em}
\begin{axiom}
  \label{axiom1}
  The invariant $\s$ is additive.  More precisely, for $6$--dimensional closed $e$-manifolds $\a$ and $\a'$,
  \begin{align*}
	\s(-\a) &= -\s(\a),\\
	\s(\a \amalg \a') &= \s(\a) + \s(\a').
  \end{align*}
\end{axiom}
\begin{axiom}
  \label{axiom2}
  If a $6$--dimensional closed $e$-manifold $\a$ bounds a $7$--dimensional $e$-mani\-fold $(Z,X,e)$,  namely
  $\bd{(Z,X,e)} \cong \a$, then
  \begin{equation*}
	\s(\a) = \Sign X.
  \end{equation*}
\end{axiom}
%\vspace{2em}
Here, $\Sign X \in \Z$ is the signature of a $4$--manifold $X$.  
The second purpose is to show
that $\s$ can detect both Haefliger's invariant (\fullref{thm:Hf}) and Milnor's triple linking
number (\fullref{thm:milnor}).

Our main results are the following.
%-------------------------------------------------------------------------------------------
\subsection{Existence and uniqueness theorem}
\label{sec:sigma}
%-------------------------------------------------------------------------------------------
All manifolds are assumed to be smooth, oriented, and compact unless otherwise stated. 
The following fundamental theorem on $6$--dimensional $e$-manifolds is the key to proving the
existence and the uniqueness of the invariant $\s$, and the \hyperlink{proof:cob}{proof} will be
given in \fullref{sec:cobordism-e}.
\begin{mainthm}\label{thm:cob}
  Any $6$--dimensional closed $e$-manifold is rationally null--cobordant.
\end{mainthm}
The statement means that, for any $6$--dimensional closed $e$-manifold $\a$, there exists a
positive integer $m$ and a $7$--dimensional $e$-manifold $\b$ such that $\bd{\b} \cong \amalg^{m} \a$ (union
of $m$ disjoint copies of $\a$). 
This is a direct consequence of the fact that the cobordism group $\O_{6}^{e}$ of
$6$--dimensional $e$-manifolds (\fullref{sec:cobordism-e}) is isomorphic to
$(\Q/\Z){}^{\oplus 2}$ (\fullref{thm:O}), and $m$ can be chosen to be the order of
the cobordism class $[\a] \in \O_{6}^{e}$ of $\a$.

The following is the existence and the uniqueness theorem of the invariant $\s$.
\begin{mainthm}\label{thm:main1}
  Let $\a$ be a $6$--dimensional closed $e$-manifold.
  Take any $7$--dimensional $e$-manifold
  $\b = (Z,X,e)$ such that $\bd{\b} \cong \amalg^{m}\a$ for some positive integer $m$
  \textup{(}such $\b$ and $m$ exist by \textup{\fullref{thm:cob}}\textup{)}.
  Then, the rational number
  \begin{equation*}
	\s(\a) \defeq \frac{\Sign X}{m}
  \end{equation*}
  depends only on the isomorphism class of $\a$.
  Moreover, the invariant $\s$ has the following properties\textup{:}
  \begin{enumerate}
	\item \label{item:axiom} The invariant $\s$ satisfies $\textup{\fullref{axiom1}}$ and
	  $\textup{\fullref{axiom2}}$.
	\item \label{item:unique} The invariant $\s$ is unique.  That is,
	  if an invariant $\s'$ of $6$--dimensional closed
	  $e$-manifolds satisfies $\textup{\fullref{axiom1}}$ and $\textup{\fullref{axiom2}}$, then
	  $\s' = \s$.
	\item \label{item:sl} If a $6$--dimensional closed $e$-manifold $\a$ bounds a $7$--dimensional
	  quasi $e$-manifold $\b = (Z,X,e)$, then
	  \[
	  \s(\a) = \Sign X- 4\L(\b).
	  \]
  \end{enumerate}
\end{mainthm}
Here, $\L(\b) \in \Q$ is the self--linking number of a $7$--dimensional quasi $e$-manifold $\b$,
and it will be defined in \fullref{sec:self-linking}.
The proof of \fullref{thm:main1} will be given in \fullref{sec:Proof-2}.
\begin{remark}
  If $\b$ is a $7$--dimensional $e$-manifold, then the formulas in \fullref{axiom2} and
  \fullref{thm:main1}\,\eqref{item:sl} are the same,  since $\L(\b) = 0$ by the definition of $\L$.
\end{remark}
\begin{remark}\label{rem:easier}
  If we use only the {axioms} to compute $\s(\a)$,
  we need to find (or construct) an $e$-manifold $\b$ such that $\bd{\b} \cong \amalg^{m} \a$
  and that the signature of the submanifold is computable, but that may not always be easy.
  However, sometimes finding a simple quasi $e$-manifold bounded by $\a$ may be much easier.
  In such cases, the formula in \fullref{thm:main1}\,\eqref{item:sl} gives us an
  alternative and effective way to compute
  $\s(\a)$.
  In fact, this formula will be used when we explorer the relationship between
  our invariant and Haefliger's invariant
  (\fullref{sec:3knot}), or Milnor's triple linking number (\fullref{sec:milnor}).
\end{remark}
The essential reason why the rational number $\s(\a)$ is independent of the choices of $\b$ is
that if a $7$--dimensional $e$-manifold $\b = (Z,X,e)$ is closed, then $\Sign X=0$
(\fullref{cor:O4-X}), and that the signature is additive with respect to the decompositions of
closed manifolds (Novikov additivity).

More generally, if $\b$ is a closed $7$--dimensional quasi $e$-manifold, then the equality
$\Sign X = 4 \L(\b)$ holds (\fullref{prop:SL}), and $\L(\b)$  is also additive with respect to the
decompositions of closed quasi $e$-manifolds (see the proof of \fullref{prop:s}).
These are the main reasons why \fullref{thm:main1}\,\eqref{item:sl} holds.

%-------------------------------------------------------------------------------------------
\subsection{An invariant of smooth embeddings}
\label{sec:emb}
%-------------------------------------------------------------------------------------------
Two manifold pairs $(Z,X)$ and $(Z',X')$ are \emph{isomorphic}
if there exists an orientation preserving diffeomorphism $f \co Z \to Z'$ such that $f(X) = X'$ as oriented
submanifolds.  The rational number $\s(W,V,e)$ defined in
\fullref{thm:main1} is not an invariant of the isomorphism class of $(W,V)$ in general, since it may depend on
the choice of $e$.
However, if we put all $e$-classes together, we obtain an invariant of $(W,V)$ as follows.
Let
\begin{equation*}
  \calE_{W,V} =
  \left\{ \left. e \in H^{2}(W \setminus V;\Q) \ \right| \
  \text{$e$ is an $e$-class of $(W,V)$} \right\}
\end{equation*}
be the set of all $e$-classes of $(W,V)$.  
For example,
if $V$ is empty, then $\calE_{W,\es} = H^{2}(W;\Q)$ by definition, and
if $V \neq \es$, then $\calE_{W,V}$ is empty or an affine subspace of $H^{2}(W
\setminus V;\Q)$ which misses the origin. 
Let
\begin{equation*}
  \s_{W,V} \co \calE_{W,V} \to \Q
\end{equation*}
be the function defined by $\s_{W,V}(e) = \s(W,V,e)$ for $e \in \calE_{W,V}$.
The following is a corollary of \fullref{thm:main1}.
\begin{maincor}\label{cor:ftn}
  For a pair $(W,V)$ of closed manifolds of dimensions $6$ and $3$, the
  function $\s_{W,V} \co \calE_{W,V} \to \Q$ is an invariant of the
  isomorphism class of $(W,V)$.
\end{maincor}
The statement means that if there is an isomorphism
$f \co (W',V') \to (W,V)$ of pair of manifolds,  then the pull--back
$f^* \co H^{2}(W \setminus V;\Q) \to H^{2}(W' \setminus V';\Q)$
restricts to a bijection
$f^* \co \calE_{W,V} \to \calE_{W',V'}$, and the identity
\[
\s_{W',V'}(f^*e) = \s_{W,V}(e)
\]
holds for any $e \in \calE_{W,V}$.

In a special case, we can obtain a rational--valued invariant of the isomorphism class of 
$(W,V)$, rather than a function--valued invariant, as follows.
\begin{definition}\label{def:simple}
  A pair $(Z,X)$ of manifolds of codimension $3$ is \emph{simple} if
  it admits at least one $e$-class and the restriction
  $H^{2}(Z;\Q) \to H^{2}(X;\Q)$ is injective.
\end{definition}
\begin{mainthm}\label{thm:Q}
  If a pair $(W,V)$ of closed manifolds of dimensions $6$ and $3$ is simple,
  then the rational number
  \begin{equation*}
	\s(W,V) \defeq \s_{W,V}(e) = \s(W,V,e),\quad e \in \calE_{W,V}.
  \end{equation*}
  is an invariant of the isomorphism class of $(W,V)$.
\end{mainthm}
The \hyperlink{proof:Q}{proof} will be given at the end of \fullref{sec:e}, and is easy.
The essential part is that $(W,V)$ is simple if, and only if, $(W,V)$ admits just one $e$-class
(\fullref{prop:eclass}\,\eqref{item:inj}).
%-------------------------------------------------------------------------------------------
\subsection{Haefliger's invariant}
\label{sec:haefliger}
%-------------------------------------------------------------------------------------------
Let $H \co \Emb(S^{3},S^{6}) \to \Z$ be Haefliger's isomorphism. 
A short review of the definition of $H$ will be given in \fullref{sec:review-h}.
Let $f \co S^{3} \to S^{6}$ be a smooth embedding, and write $M_{f} = f(S^{3})$.

There is an easy--to--check condition for the simplicity of pairs of manifolds as follows.
\begin{proposition}\label{prop:simple}
  Let $(Z,X)$ be a pair of manifolds of codimension $3$, and
  assume that the restriction $H^{2}(Z;\Q) \to H^{2}(X;\Q)$ is an
  isomorphism.  Then, $(Z,X)$ is simple if, and only if, $(\pair{X})$ is rationally
  null--homologous in $(\pair{Z})$.
\end{proposition}
The proof will be given in \fullref{sec:e}.
By \fullref{prop:simple}, the pair $(S^{6},M_{f})$ is simple,
and the rational number $\s(S^{6},M_{f})$ is well--defined by \fullref{thm:Q}.
The relationship between Haefliger's invariant $H(f)$ and our invariant $\s(S^{6},M_{f})$ is the
following
one.
\begin{mainthm}\label{thm:Hf}
  For a smooth embedding $f \co S^{3} \to  S^{6}$, we have
  \begin{equation*}
	\s(S^{6},M_{f}) =  -8 H(f).
  \end{equation*}
\end{mainthm}
The proof will be given in \fullref{sec:proof5}, and
it will turn out that our invariant $\s$ is a natural generalization of Haefliger's invariant $H$.

There are some generalizations of Haefliger's invariant due to
Takase \cite{takase-hom3} and \linebreak Skopenkov \cite{skopenkov}.
Takase~\cite{takase-hom3} \cite{takase-geom} proved that there is a bijection
$\O \co \Emb(M,S^{6}) \to \Z$ for any integral homology $3$--sphere $M$ such that if $M=S^{3}$ then
$\O=H$ .
Our invariant recovers Takase's invariant too (\fullref{cor:takase}), and that is a direct
consequence of the geometric formula for $\s(W,V,e)$ (\fullref{thm:geometric-formula}):
\begin{equation*}
  \s(W,V,e) = 
  \Sign S - \int_{S}^{}e(\nu_{S})^{2}
\end{equation*}
Here, $S \supset W$ is a Seifert surface of $V$ which Poincar\'e dual is $e/2$, and which rational normal
Euler class $e(\nu_{S}) \in H^{2}(S;\Q)$ is trivial over $\bd{S}$. 
If $W=S^{6}$ and $V$ is an integral homology $3$--sphere, then the right--hand side is nothing but ($-8$
times) the definition of $\O$. 

Recently, Skopenkov~\cite{skopenkov} proved a classification theorem of elements in
$\Emb(M,S^{6})$ for any oriented connected closed $3$--manifolds $M$.  
When $M = S^{3}$, his invariant 
$\mu \co \operatorname{Wh}^{-1}(0) \to \Z$ (where
$\operatorname{Wh}^{-1}(0)$ is a subset\footnote{The subset $\operatorname{Wh}^{-1}(0) \subset
\Emb(M,S^{3})$ consists of  
elements which Whitney invariant $\operatorname{Wh} \co \Emb(M,S^{6}) \to H_{1}(M;\Z)$
\cite{skopenkov} vanish.}
of $\Emb(M,S^{6})$, and $\mu$ is called the Kreck
invariant in his paper) coincides with $H$.
The invariants $\mu$ and $\s$ seem to be closely related, and possibly identical
(up to multiplication by a constant) for any $M$.
%-------------------------------------------------------------------------------------------
\subsection{Milnor's triple linking number}
\label{sec:triple}
%-------------------------------------------------------------------------------------------
Let $L$ be an oriented algebraically split $3$--component link in $\R^{3}$, and let
$(T^{3}\x S^{3}, M_{L})$ be the manifold pair defined as before.
In \fullref{sec:milnor}, we will prove that $(T^{3}\x S^{3}, M_{L})$ is simple (\fullref{prop:link-e}),
and consequently, the rational number $\s(T^{3}\x S^{3}, M_{L})$ is well--defined by \fullref{thm:Q}.
\begin{remark}\label{rem:link-homotopy}
  It is easy to see that $\s(T^{3}\x S^{3}, M_{L})$ is a
  link homotopy invariant of $L$ without explicit computations, in fact, we can see that the
  isotopy class of the submanifold $M_{L}$ depends only on the link homotopy type of
  $L$ as follows.
  Suppose that two
  algebraically split $3$--component links $L$ and $L'$ in $\R^{3}$ have the same
  link homotopy type, and
  let $\left\{ L(t) \right\}_{t \in [0,1]}$ be a smooth link homotopy\footnote{Each connected
  component $K_{i}(t)$ of
  each intermediate link $L(t) = K_{1}(t) \cup K_{2}(t) \cup K_{3}(t)$ may intersects itself, but no other
  components $K_{j}(t)$ ($i \neq j$).}
  from $L$ to $L'$.
  For each $t \in [0,1]$, we can construct a smoothly embedded $3$--submanifold
  $M_{L(t)} \subset T^{3} \x S^{3}$, in exactly the same way as we construct $M_{L}$. 
  The obtained family $\left\{ M_{L(t)} \right\}_{t \in [0,1]}$ is a
  smooth isotopy from $M_{L}$ to $M_{L'}$.
\end{remark}
The relationship between Milnor's triple linking number $\mu(L)$ and our invariant\linebreak
$\s(T^{3}\x S^{3}, M_{L})$ is the following one.
\begin{mainthm}\label{thm:milnor}
  For an oriented algebraically split $3$--component link $L$ in $S^{3}$, we have
  \[
  \s(T^{3} \x S^{3}, M_{L}) = -8 \mu(L).
  \]
\end{mainthm}
The proof will be given in \fullref{sec:milnor}.

%\vspace{1.5em}

Now, here is the plan of the paper.
In \fullref{sec:preliminaries},
we introduce definitions and notation which are
necessary to understand the main theorems given in this section.
In \fullref{sec:BSO-KQ2},
we study some elementary facts on the low--dimensional oriented cobordism groups
$\O_{*}(K(\Q,2))$ and
$\O_{*}(BSO(3))$ of the Eilenberg--MacLane space
$K(\Q,2)$ of type $(\Q,2)$ and the classifying space $BSO(3)$ of the Lie group $SO(3)$, 
and this is a preliminary to the next section.
In \fullref{sec:cobordism-e}, we show that there is an isomorphism $\O_{6}^{e} \cong (\Q/\Z)^{\oplus 2}$
(\fullref{thm:O}), and we prove \fullref{thm:cob} as a consequence of this isomorphism.
The isomorphism is given by a short exact sequence
\[
0 \to \O_{4}(BSO(3)) \to \O_{6}(K(\Q,2)) \to \O_{6}^{e} \to 0,
\]
which is isomorphic to $0 \to \Z^{\oplus 2} \to \Q^{\oplus 2} \to (\Q/\Z)^{\oplus 2} \to
0$ (see the proof of \fullref{thm:O}).
\fullref{sec:proof1} is devoted to the proof of \fullref{thm:main1}, roughly speaking, which
relies on the two properties of $e$-manifolds as follows:
\begin{enumerate}
  \item  \fullref{thm:cob} implies the existence and the uniqueness of the invariant $\s$.
  \item  The formula $\Sign X = \L(\b)$ (\fullref{prop:SL}) implies that $\s$ is well--defined,
	and that \fullref{thm:main1}\,\eqref{item:sl} holds.
\end{enumerate}
\fullref{sec:e} is the study of necessary and sufficient conditions ensuring the existence and
uniqueness of $e$-classes.
In particular, we prove that a pair $(W,V)$ is simple if, and only if, it is admits just one $e$-class
(\fullref{prop:eclass}\,\eqref{item:inj}).
The proof of \fullref{thm:Q} is given at the end of the section.
In \fullref{sec:seifert-intro}, we study the relationship between Seifert surfaces and $e$-classes.
In \fullref{sec:3knot}, we prove \fullref{thm:Hf},
and it will turn out that our invariant $\s$ is a natural generalization of the Haefliger's invariant $H$.
We also prove the geometric formula (\fullref{thm:geometric-formula}) for $\s(W,V,e)$ when $e/2$ is represented by a Seifert
surface of $V$, and as a direct consequence, we prove that our invariant also recovers Takase's invariant
(\fullref{cor:takase}).
In \fullref{sec:milnor}, we prove \fullref{thm:milnor}.

%%%%%%%%%%%%%%%%%%%%%%%%%%%%%%%%%%%%%%%%%%%%%%%%%%%%%%%%%%%%%%%%%%%%%%%%%%%%%%%
%%%   Acknowledgments
\textit{Acknowledgments.}
I would like to thank Professor
Mikio~Furuta,
Toshitake~Kohno,
and
Christine~Lescop
for their
advice and support.
I also would like to thank
Osamu~Saeki,
Masamichi~Takase,
and
Tadayuki~Watanabe
for helpful suggestions and comments, and
Michael~Eisermann
for revising the English text.
%%%%%%%%%%%%%%%%%%%%%%%%%%%%%%%%%%%%%%%%%%%%%%%%%%%%%%%%%%%%%%%%%%%%%%%%%%%%%%%

%%%%%%%%%%%%%%%%%%%%%%%%%%%%%%%%%%%%%%%%%%%%%%%%%%%%%%%%%%%%%%%%%%%%%%%%%%%%%%%%%%%%%%%%%
\section{Preliminaries}
\label{sec:preliminaries}
%%%%%%%%%%%%%%%%%%%%%%%%%%%%%%%%%%%%%%%%%%%%%%%%%%%%%%%%%%%%%%%%%%%%%%%%%%%%%%%%%%%%%%%%%
%-------------------------------------------------------------------------------------------
\subsection{Notation}
\label{sec:notation}
%-------------------------------------------------------------------------------------------
We use the ``outward normal first'' convention for boundary orientation of manifolds.
For an oriented real vector bundle $E$ of rank $3$ over a manifold
$X$, we denote the associated unit sphere bundle by
$\rho_{E} \co S(E) \to X$,
and let $F_{E} \subset TS(E)$ denote the vertical tangent subbundle of $S(E)$ with respect to $\rho_{E}$.
The orientations of $F_{E}$ and $S(E)$ are given by the isomorphisms
$\rho_{E}^*E \cong \R_{E} \oplus F_{E}$ and $TS(E) \cong \rho_{E}^*TX \oplus
F_{E}$,
where $\R_{E} \subset \rho_{E}^*E$ is the tautological real line bundle of $E$ over $S(E)$.
Consequently, the Euler class
\[
e(F_{E}) \in H^{2}(S(E);\Z)
\]
of $F_{E}$ is defined.

Next, let $(Z,X)$, $Z \supset X$, be a pair of manifolds, and
we assume that $X$ is properly embedded in $Z$ and the codimension is $3$.
Throughout this paper, we always impose these assumptions for all pairs of manifolds.
In particular, when we write
$(W,V)$ or $(Z,X)$, we always mean a pair of manifolds of codimension $3$
such that $W$ and $V$ are closed and $\dim W = 6$,
and that $Z$ and $X$ may have boundaries and $\dim Z$ can be any (mainly assumed to be $7$ or
$6$).
Denote by $\nu_{X}$ the normal bundle of $X$, which can be identified with a tubular
neighborhood of $X$ so that $X \subset \nu_{X} \subset Z$.
For simplicity, we will write
\begin{equation*}
  \HX = S(\nu_{X}), \qquad
  \rho_{X} = \rho_{\nu_{X}} \co \HX \to X, \qquad
  F_{X} = F_{\nu_{X}}.
\end{equation*}

Let us write $(W,V) = \bd{(Z,X)}$ for the boundary pair of $(Z,X)$ for a moment.
We can define $\nu_{V}$, $F_{V}$, $\HV$, $\rho_{V} \co \HV \to V$,
etc.~in exactly the same way as above.
Let $\left( (0,1] \x W, (0,1] \x V\right)$ be the pair of collar neighborhoods of the boundary pair
$\left( \left\{ 1 \right\} \x W, \left\{ 1 \right\} \x V \right) = (W,V)$. 
Without loss of generality, we shall always assume
$\nu_X|_{(0,1] \x V} = (0,1] \x \nu_{V}$
as tubular neighborhoods of $(0,1] \x V$ in $(0,1] \x W$.
Consequently, we have
\[
	\bd{\HX} = \HV,\qquad
	e(F_{X})|_{\HV} = e(F_{V}).
\]

%-------------------------------------------------------------------------------------------
\subsection{\texorpdfstring{$e$-classes and $e$-manifolds}{e-classes and e-manifolds}}
\label{sec:emfd}
%-------------------------------------------------------------------------------------------
Here is the definition of $e$-class and quasi $e$-class.
\begin{definition}\label{def:eclass}
  Let $(Z,X)$ be a manifold pair of dimensions $n$ and $n-3$.
  \begin{enumerate}
	\item A cohomology class $e \in H^{2}(Z \setminus X;\Q)$ is called an \emph{$e$-class} of
	  $(Z,X)$ if
	  \begin{enumerate}
		\item $e|_{\HX} = e(F_{X})$ over $\Q$. 
	  \end{enumerate}
	  We call $(Z,X,e)$ an $n$--dimensional \textup{$e$-manifold}.
	\item\label{item:qe} A cohomology class $e \in H^{2}(Z \setminus X;\Q)$ is called a \emph{quasi $e$-class}
	  of $(Z,X)$ if 
	  \begin{enumerate}
		\item $e|_{\bd{Z} \setminus \bd{X}}$ is an $e$-class of $\bd{(Z,X)}$, and
		\item $\bk{[S^{2}_{p}], e} = 2$ for all $p \in X$. 
	  \end{enumerate}
	  We call $(Z,X,e)$ an $n$--dimensional \textup{quasi $e$-manifold}.
  \end{enumerate}
\end{definition}
Here, $S^{2}_{p} = \rho_{X}^{-1}(p) \subset \HX$ is the fiber of $\rho_{X}$ at $p$, and the
bracket $\bk{~,~}$ denotes the pairing
of a homology class and a cohomology class.
Note that any $e$-class is a quasi $e$-class,
since $\bk{[S^{2}_{p}], e(F_{X})} = 2$, which is the Euler characteristic of the $2$--sphere.
Also note that the boundary of a quasi $e$-manifold is an $e$-manifold by definition.

For (quasi) $e$-manifolds $\b = (Z,X,e)$ and $\b' = (Z',X',e')$,
if there exists an isomorphism $f \co (Z',X') \to (Z,X)$ of
pair of manifolds such that $f^*e
=e'$, then we say $\b$ and $\b'$ are \emph{isomorphic} (denoted by $\b \cong \b'$).
The empty $e$-manifold $(\es,\es,0)$, where $0 \in H^{2}(\es \setminus \es;\Q)$,
will be simply denoted by $\es$.
If $\bd{\b} \cong \es$, then we say $\b$ is \emph{closed}.
If a closed $e$-manifold $\a$ bounds an $e$-manifold $\b$, i.e.~$\bd{\b} \cong
\a$, then we say $\a$ is \emph{null--cobordant}.
If $\a \amalg (-\a')$ is null--cobordant, then we say $\a$ and $\a'$ are
\emph{cobordant}.

%-------------------------------------------------------------------------------------------
\subsection{\texorpdfstring{Self--linking form and self--linking number}{Self-linking form and
self-linking number}}
\label{sec:self-linking}
%-------------------------------------------------------------------------------------------
For a quasi $e$-manifold $\b = (Z,X,e)$,
we define the \emph{self--linking form} $\g \in H^{2}(X;\Q)$
and the \emph{self--linking number} $\L(\b) \in \Q$ as follows.
Here, $\L(\b)$ is defined only when $\dim \b = 7$, i.e.~$\dim X = 4$.

By the Thom--Gysin exact sequence
\begin{equation*}
	\dotsb \to 0 \to H^{2}(X;\Q) \xrightarrow{\rho_{X}^*} H^{2}(\HX;\Q)
	\xrightarrow{\rho_{X}{}_{!}} H^{0}(X;\Q) \to \dotsb
\end{equation*}
of $\nu_{X}$,
the cohomology class $e|_{\HX} - e(F_{X}) \in H^{2}(\HX;\Q)$ belongs to the image of the
pull--back $\rho_{X}^*$ which is injective.
We define the self--linking
form $\g \in H^{2}(X;\Q)$ of $\b$ to be the unique cohomology class such that
\begin{equation*}
  e|_{\HX} = e(F_{X}) +  2\, \rho_{X}^* \g.
\end{equation*}

Since $\bd{\b}$ is an $e$-manifold, we have $\g|_{\bd{X}} = 0$.
There exists $\Tg \in H^{2}(\pair{X};\Q)$ such that the homomorphism
$H^{2}(\pair{X};\Q) \to H^{2}(X;\Q)$ maps $\Tg$ to $\g$.
When $\dim \b = 7$, the self--linking number $\L(\b)$ of $\b$ is defined by
\begin{equation*}
  \L(\b) = \int_{X}\Tg^{2}.
\end{equation*}
It is easy to check that $\L(\b)$ does not depend on the choice of $\Tg$.
Note that $\g = 0$ if $\b$ is an $e$--manifold, so $\L(\b) = 0$ if $\b$ is a
$7$--dimensional $e$-manifold.

  In \fullref{prop:seifert},
we will give an interpretation of $\g$ by using a
Seifert surface of $X$.

%%%%%%%%%%%%%%%%%%%%%%%%%%%%%%%%%%%%%%%%%%%%%%%%%%%%%%%%%%%%%%%%%%%%%%%%%%%%%%%%%%%%%%%%%%
\section{\texorpdfstring{Oriented cobordism groups of $BSO(3)$ and $K(\Q,2)$}{Oriented
cobordism groups of BSO(3) and K(Q,2)}}
\label{sec:BSO-KQ2}
%%%%%%%%%%%%%%%%%%%%%%%%%%%%%%%%%%%%%%%%%%%%%%%%%%%%%%%%%%%%%%%%%%%%%%%%%%%%%%%%%%%%%%%%%%
Let $K(\Q,2)$ be the Eilenberg--MacLane space of type $(\Q,2)$,
i.e.~$\pi_{2}(K(\Q,2)) \cong \Q$ and $\pi_{i}(K(\Q,2)) = 0$ for $i \neq 2$,
and $BSO(3)$ the classifying space of the Lie group $SO(3)$.
We can assume that $BSO(3)$ and $K(\Q,2)$ have structures of CW--complexes.
Let $\O_{*}(Y)$ denote the oriented cobordism group of a CW--complex $Y$.
As a preparation for the next section, in this section
we study some relationship between $\O_{*}(BSO(3))$ and $\O_{*}(K(\Q,2))$.
%-------------------------------------------------------------------------------------------
\subsection{Homology groups}
\label{sec:homology}
%-------------------------------------------------------------------------------------------
We begin by recalling some elementary facts on the homology groups of $K(\Q,2)$ and
$BSO(3)$.
The homotopy class of a map $\CP^{\infty} \to  K(\Q,2)$ from the infinite dimensional complex
projective space $\CP^{\infty}$ ($\simeq K(\Z,2)$),
corresponding to the inclusion $\Z \hookrightarrow \Q$ on
the second homotopy groups, provides an isomorphism
between the reduced homology groups (cf.~\cite{GM}):
%(cf. \cite[Theorem\,7.7]{GM}).
\begin{equation}\label{eq:CP}
  \begin{split}
	\TH_{k}(K(\Q,2);\Z) &\cong \TH_{k}(\CP^{\infty};\Q)\\
	&\cong
	\begin{cases}
	  \Q &\text{if $k$ is positive even}\\
	  0 &\text{otherwise}
	\end{cases}
  \end{split}
\end{equation}
Let $a_{1} \in H^{2}(K(\Q,2);\Q) \cong \Q$ denote the element dual to
$1 \in \pi_{2}(K(\Q,2)) \cong H_{2}(K(\Q,2);\Q)$, then
the $k$--th power $a_{1}^{k}$ of $a_{1}$ generates
$H^{2k}(K(\Q,2);\Q)$ over $\Q$ for $k \geq 1$.

It is easy to check that the low--dimensional homology groups of $BSO(3)$ are given as follows:
\begin{equation}\label{tbl:BG}
  \begin{array}{c|ccccccc}
	k & 0 & 1 & 2 & 3 & 4 & 5\\
	\hline
	H_{k}(BSO(3);\Z) & \Z & 0 & \Z/2 & 0 & \Z & \Z/2
  \end{array}
\end{equation}
This table, for example, is obtained by use of the Serre spectral sequence of the universal
principal $SO(3)$-bundle, and by the fact that the cohomology ring $H^{*}(BSO(3);\Z/2)$
is a free polynomial algebra generated by the second and third Stiefel-Whitney
classes over $\Z/2$.
The group $H_{4}(BSO(3);\Z)$ is generated by the dual element of the first Pontryagin
class $p_{1} \in H^{4}(BSO(3);\Z)$.

%-------------------------------------------------------------------------------------------
\subsection{Cobordism groups}
\label{sec:oriented-cobordism}
%-------------------------------------------------------------------------------------------
The next step is the study of the low--dimensional oriented cobordism groups of $K(\Q,2)$
and $BSO(3)$.
In low--dimensions, the cobordism group $\O_* = \O_{*}(pt)$ of one point $pt$ is given as
follows (cf.~\cite[Section\,17]{Milnor-characteristic}):
\begin{equation}\label{tbl:O-pt}
  \begin{array}{c|ccccccc}
	k & 0 & 1 & 2 & 3 & 4 & 5 & 6 \\
	\hline
	\O_{k} & \Z & 0 & 0 & 0 & \Z & \Z/2 & 0
  \end{array}
\end{equation}
Here, the isomorphism $\O_{4} \cong \Z$ is given by the signature of $4$--manifolds.

In general, for any CW--complex $Y$, the Atiyah--Hirzebruch spectral sequence
$E^{n}_{p,q}(Y)$ for $\O_{*}(Y)$ converges
(cf.~\cite[Theorem\,15.7]{switzer}):
\[
E^{2}_{p,q}(Y) = H_{p}(Y;\O_{q})  \ \Longrightarrow \ \O_{p+q}(Y)
\]
The following lemma is an easy application of the Atiyah--Hirzebruch spectral sequence.
\begin{lemma}\label{lem:O}
  The following isomorphisms hold\textup{:}
  \begin{equation*}
	\O_{6}(K(\Q,2)) \cong \Q^{\oplus 2},\qquad
	\O_{3}(BSO(3)) = 0,\qquad
	\O_{4}(BSO(3)) \cong \Z^{\oplus 2}.
  \end{equation*}
\end{lemma}
\begin{proof}
  We use the isomorphism \eqref{eq:CP} and the tables  \eqref{tbl:BG} and \eqref{tbl:O-pt} to 
  prove this lemma.
  The Atiyah--Hirzebruch spectral sequence $E^{n}_{p,q} = E^{n}_{p,q}(K(\Q,2))$ converges on
  the $E^{2}$-stage within the range $p+q \leq 6$, and so $E^{\infty}_{p,q} \cong E^{2}_{p,q}$
  in the same range.
  Consequently, we have
  \begin{equation*}
	E^{\infty}_{p,6-p} \cong
	\begin{cases}
	  \Q & \text{if $p=6,2$},\\
	  0 & \text{otherwise},
	\end{cases}
  \end{equation*}
  and therefore, $\O_{6}(K(\Q,2)) \cong \Q^{\oplus 2}$.

  Similarly, the spectral sequence $F^{n}_{p,q} = E^{n}_{p,q}(BSO(3))$
  converges on the $F^{2}$-stage in the range $p+q \leq 4$, and 
  \begin{equation*}
	F^{\infty}_{p,4-p} \cong
	\begin{cases}
	  \Z & \text{if $p=4,0$,}\\
	  0 & \text{otherwise.}
	\end{cases}
  \end{equation*}
  Thus, $\O_{4}(BSO(3)) \cong \Z^{\oplus 2}$.
  The vanishing of $\O_{3}(BSO(3))$ follows from $F^{2}_{p,3-p} = 0$ for all $p$.
\end{proof}
A pair $(W,e)$ of a closed $6$--manifold $W$ and a cohomology class $e \in H^{2}(W;\Q)$
represents a cobordism
class $[W,e] \in \O_{6}(K(\Q,2))$.  Here, we identify $e$ with the homotopy class of a map $f \co W \to
K(\Q,2)$ such that $f^*a_{1} = e$.
Define a homomorphism
$\chi \co \O_{6}(K(\Q,2)) \to \Q^{\oplus 2} $
by
\begin{align*}
  \chi([W,e]) &= \left(\chi_{1}(W,e), \chi_{2}(W,e) \right),\\
  \chi_{1}(W,e) &= \frac{1}{6}\int_{W}p_{1}(TW)\,e - e^{3} \in \Q,\\
  \chi_{2}(W,e) &= \frac{1}{2} \int_{W}e^{3} \in \Q.
\end{align*}
Similarly, a pair $(X,E)$ of a closed $4$--manifold $X$ and an oriented vector bundle $E$ of rank $3$
over $X$ represents a cobordism class $[X,E] \in \O_{4}(BSO(3))$.  Here, we identify the
isomorphism class of $E$ with
the homotopy class of the classifying map $X \to BSO(3)$ of $E$.
Define a homomorphism $\xi \co \O_{4}(BSO(3)) \to \Z^{\oplus 2}$ by
\begin{equation*}\label{eq:def-xi}
  \xi([X,E]) = \left(  \Sign X, \int_{X}p_{1}(E))\right).
\end{equation*}
We will see soon that the homomorphisms $\chi$ and $\xi$ are isomorphic (\fullref{lem:chi-xi}).

%-------------------------------------------------------------------------------------------
\subsection{\texorpdfstring{Homomorphism $\O_{4}(BSO(3)) \to \O_{6}(K(\Q,2))$}{Homomorphism
O4(BSO(3)) --> O6(K(Q,2))}} 
\label{sec:hom-O}
%-------------------------------------------------------------------------------------------
Let us consider the homomorphism
\[
\up \co \O_{4}(BSO(3)) \to \O_{6}(K(\Q,2))
\]
defined by
$
\up([X,E]) = [S(E), e(F_{E})]
$
for $[X,E] \in \O_{4}(BSO(3))$.
For a pair $(X,E)$ representing an element in $\O_{4}(BSO(3))$,
the characteristic classes of the vector bundles $E$, $F_{E}$, $TX$, and $TS(E)$
satisfy the following relations:
\begin{align}
  \label{eq:pe1}
  e(F_{E})^{2} &= p_{1}(F_{E}) = \rho_{E}^*p_{1}(E)\\
  \label{eq:pe2}
  &\equiv p_{1}(TS(E)) - \rho_{E}^*p_{1}(TX)\quad \text{(modulo $2$--torsion
  elements)},\\
  \label{eq:gysin}
  {\rho_{E}}_{!} e(F_{E}) &=2
\end{align}
Here, ${\rho_{E}}_{!} \co H^{2}(S(E);\Z) \to H^{0}(X;\Z)$ is the
Gysin homomorphism of $\rho_{E}$, and  $2 \in
H^{0}(X;\Z)$ denotes the element given by the constant function on $X$ with
the value $2$.
The Hirzebruch signature theorem states that
\begin{equation}
  \label{eq:hirzebruch}
  \Sign X = \frac{1}{3}\int_{X}p_{1}(TX).
\end{equation}
The next two lemmas are easy to prove.
\begin{lemma}
  \label{lem:cu}
  $\chi \up = \xi$.  Namely, for any pair $(X,E)$ of closed $4$--manifold $X$ and an oriented
  vector bundle $E$ of rank $3$ over $X$, we have
  \begin{equation*}
	\label{eq:SO}
	\chi([S(E), e(F_{E})]) = \left( \Sign X,\  \int_{X}p_{1}(E) \right).
  \end{equation*}
\end{lemma}
\begin{proof}
  This follows from the formulas \eqref{eq:pe1}, \eqref{eq:pe2},
  \eqref{eq:gysin}, and \eqref{eq:hirzebruch}.
  In fact, these imply
  \begin{equation*}
	p_{1}(TS(E)) e(F_{E}) - e(F_{E})^{3} = \rho_{E}^*p_{1}(TX)\, e(F_{E})
  \end{equation*}
  over $\Q$, and
  \begin{equation*}
	\chi_{1}(S(E), e(F_{E})) =
	\frac{1}{6}\int_{S(E)}^{}\rho_{E}^*p_{1}(TX)\, e(F_{E}) =
	\frac{1}{3}\int_{X}^{}p_{1}(TX) =
	\Sign X.
  \end{equation*}
  The equality $\chi_{2}(S(E),e(F_{E})) = \int_{X}^{}p_{1}(E)$ can be obtained in a similar way.
\end{proof}
\begin{lemma}
  \label{lem:chi-xi}
  The homomorphisms $\chi \co \O_{6}(K(\Q,2)) \to \Q^{\oplus 2}$ and $\xi \co
  \O_{4}(BSO(3)) \to \Z^{\oplus 2}$  are isomorphisms.
\end{lemma}
\begin{proof}
  For $k=0,1$, let $F_{k}$ be an oriented vector bundle of rank $2$ over $\CP^{2}$ such
  that $\bk{[\CP^{1}], e(F_{k})} = k$, and we set $u_{k} = [\CP^{2}, F_{k} \oplus \R]
  \in \O_{4}(BSO(3))$.
  Then, two elements $\xi(u_{0}) = (1,0)$ and $\xi(u_{1}) = (1,1)$ form a basis
  of the abelian group $\Z^{\oplus 2}$.
  Therefore, $\xi$ is a surjective homomorphism from $\O_{4}(BSO(3)) \cong \Z^{\oplus 2}$
  (\fullref{lem:O}) to $\Z^{\oplus 2}$.
  This means that $\xi$ is an isomorphism.

  Similarly, we have
  $\chi(\up(u_{0})) = (1,0)$ and $\chi(\up(u_{1})) = (1,1)$ by \fullref{lem:cu},
  and these two elements form a basis of the vector space $\Q^{\oplus 2}$.
  Therefore, $\chi$ is a linear homomorphism from $\O_{6}(K(\Q),2) \cong \Q^{\oplus
  2}$ (\fullref{lem:O}) to $\Q^{\oplus 2}$ of rank $2$.
  This means that $\chi$ is an isomorphism.
\end{proof}
The following proposition is the goal of this section.
\begin{proposition}
  \label{prop:exact-O}
  The sequence of abelian groups
  \begin{equation*}
	\label{eq:exact-O}
	0 \to \O_{4}(BSO(3)) \xrightarrow{\up} \O_{6}(K(\Q,2)) \xrightarrow{\chi'}
	\left( \Q/\Z \right)^{\oplus 2} \to 0
  \end{equation*}
  is exact, where $\chi' = \chi \mod{\Z^{\oplus 2}}$.
\end{proposition}
\begin{proof}
  This follows from that, the diagram
  \begin{equation*}
	\begin{CD}
	  \O_{4}(BSO(3)) @>{\up}>> \O_{4}(K(\Q,2))\\
	  @V{\xi}V{\cong}V @V{\chi}V{\cong}V\\
	  \Z^{\oplus 2} @>\text{inclusion}>> \Q^{\oplus 2}
	\end{CD}
  \end{equation*}
  commutes (\fullref{lem:cu}) and the vertical arrows are isomorphic (\fullref{lem:chi-xi}).
\end{proof}

%%%%%%%%%%%%%%%%%%%%%%%%%%%%%%%%%%%%%%%%%%%%%%%%%%%%%%%%%%%%%%%%%%%%%%%%%%%%%%%%%%%%%%%%%%%%%
\section{\texorpdfstring{Cobordism group of $6$--dimensional
$e$-manifolds}{Cobordism group of 6-dimensional e-manifolds}}
\label{sec:cobordism-e}
%%%%%%%%%%%%%%%%%%%%%%%%%%%%%%%%%%%%%%%%%%%%%%%%%%%%%%%%%%%%%%%%%%%%%%%%%%%%%%%%%%%%%%%%%%%%%
We define $\O_{6}^{e}$ to be the cobordism group of $6$--dimensional $e$-manifolds, namely, it is
an abelian group consisting of the cobordism classes $[\a]$ of $6$--dimensional closed
$e$-manifolds $\a$, with the group structure given by the disjoint sum.  Note that $[\a] + [\a'] = [\a \amalg
\a']$, $-[\a] = [-\a]$, and $0 = [\es]$.
In this section, we prove that $\O_{6}^{e}$ is isomorphic to
$(\Q/\Z)^{\oplus 2}$ (\fullref{thm:O}), and then we prove \fullref{thm:cob}.

We begin by preparing some notation as follows.
For a pair $(Z,X)$, we will write
\begin{equation}
  \label{eq:ZX}
  Z_{X} = Z \setminus U_{X},
\end{equation}
where $U_{X}$ is the total space of the open unit disk bundle of $\nu_{X}$.
If $(W,V) = \bd{(Z,X)}$ denotes the boundary pair,
then the manifold $W_{V} = W \setminus U_{V}$ can be defined in the same way as above.
In line with our orientation conventions,
the boundaries of $Z_{X}$ and $W_{V}$ are given as follows:
\begin{equation}
  \label{eq:ori-bdry}
  \bd{Z_{X}} = W_{V} \cup (\mp \HX),\qquad
  \bd{W_{V}} = \pm \HV
\end{equation}
Here, the symbols $\mp = (-1)^{\dim Z}$ and $\pm = (-1)^{\dim W}$ are the signs of
orientations.
Note that $Z_{X}$ have the corner $\HV$ which is empty when $X$ is closed.

%-------------------------------------------------------------------------------------------
\subsection{\texorpdfstring{Extension of the cobordism group of $6$--dimensional
$e$-manifolds}{Extension of the cobordism group of 6-dimensional e-manifolds}}
\label{sec:ext}
%-------------------------------------------------------------------------------------------
In this subsection, we show that
any element in $\O_{6}^{e}$ can be represented by a $6$--dimensional closed $e$-manifold with empty
submanifold.  More precisely, let
\[
\pi \co \O_{6}(K(\Q,2)) \to \O_{6}^{e}
\]
be the homomorphism defined by $\pi([W,e]) = [W,\es,e]$ for $[W,e] \in \O_{6}(K(\Q,2))$, and
we prove that $\pi$ is surjective (\fullref{prop:pi-onto}).

For a $6$--dimensional closed $e$-manifold $\a = (W,V,e)$, we construct a cobordism class
$[W',e'] \in \O_{6}(K(\Q,2))$ such that $\pi([W',e']) = [\a]$ as follows.
Since the oriented cobordism group $\O_{3}(BSO(3))$ vanishes by \fullref{lem:O},
there exists a pair $(X,E)$ of a $4$--manifold $X$ and an
oriented vector bundle $E$ of rank $3$ over $X$ equipped with fixed identifications
$\bd{X} = V$ and $E|_{V} = \nu_{V}$.
Two pairs $(S(E),e(F_{E}))$ and $(W_{V}, e|_{W_{V}})$ have the common boundary
\begin{equation*}
  \bd{(S(E), e(F_{E}))} = (\HV, e(F_{V})) =
  \bd{(W_{V}, e|_{W_{V}})}.
\end{equation*}
Let us consider the closed $6$--manifold
\begin{equation*}
  W' = W_{V} \cup_{\HV} (-S(E))
\end{equation*}
obtained from $W_{V}$ and $-S(E)$ by gluing along the common boundaries
(namely, $W'$ is constructed by performing a kind of surgery along $V$, replacing the tubular
neighborhood of $V$ with $-S(E)$). 
There exists $e' \in H^{2}(W';\Q)$
such that $e'|_{W_{V}} = e|_{W_{V}}$ and $e'|_{S(E)} = e(F_{E})$, and we have a cobordism class
\[
[W',e'] \in \O_{6}(K(\Q,2)).
\]
\begin{proposition}
  \label{prop:pi-onto}
  We have $\pi([W',e']) = [\a]$ in $\O_{6}^{e}$.
  Consequently, the homomorphism $\pi \co \O_{6}(K(\Q,2)) \to \O_{6}^{e}$ is surjective.
\end{proposition}
\begin{proof}
  We only need to show the existence of a $7$--dimensional $e$-manifold $\b$ bounded by $\a
  \amalg (-\a')$, where $\a' = (W',\es,e')$.
  Let $I = [0,1]$ be the interval.
  In this proof, for a subset $A \subset W$, we write $A_{t} = \left\{ t \right\} \x A \subset I \x W$ for $t = 0, 1$.
  The boundaries of the $7$--manifold $I \x W$ and the closed unit disk bundle $D(E)$ of $E$ are given as
  follows:
  \begin{align*}
	\bd{(I \x W)} &=  \left( -W_{0} \right) \amalg W_{1}\\
	\bd{D(E)} &= S(E) \cup D(\nu_{V})
  \end{align*}

  Gluing the manifolds $I \x W$ and $D(E)$ along $D(\nu_{V})_{0} \subset W_{0}$ and $D(\nu_{V}) \subset
  \bd{D(E)}$ by the identity map,
  we obtain a $7$--manifold
  \begin{equation*}
	Z = D(E) \cup_{D(\nu_{V})_{0}} \left( I \x W \right)
  \end{equation*}
  with the boundary
  \begin{equation*}
	\begin{split}
	  \bd{Z} &= W_{1} \amalg \left( S(E) \cup_{\HV_{0}} (- (W_{V})_{0})  \right)\\
	  &\cong W \amalg (-W'),
	\end{split}
  \end{equation*}
  and we shall assume that $\bd{Z}$ is smooth after the corner $\HV_{0}$ is rounded,
  see \fullref{fig:gluing}.
  \begin{figure}[htbp]
	\begin{center}
	  %\input{glue0.pstex_t}
	  %%%%%%%%%%%%%%%%%%%%%%%%%%%%%%%%%%%%%%%%%%%%%%%%%%%%%%%%%%%%%%%%%%%%%%%

	  \begin{picture}(0,0)%
		\includegraphics{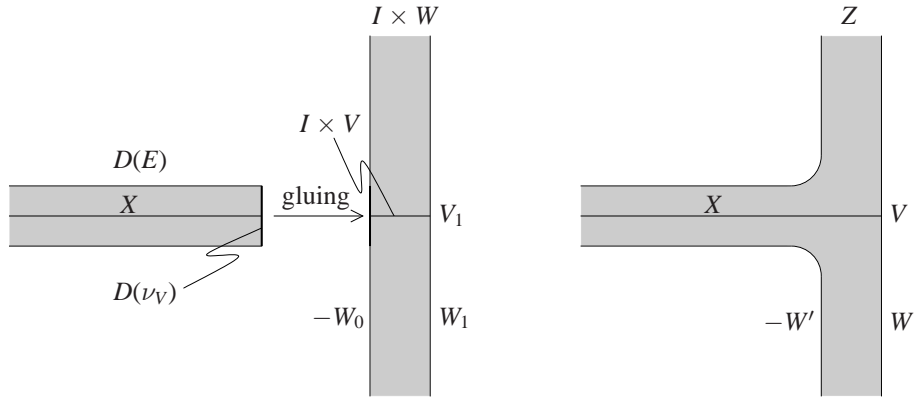}%
	  \end{picture}%
	  \setlength{\unitlength}{3315sp}%
	  \begingroup\makeatletter\ifx\SetFigFont\undefined%
	  \gdef\SetFigFont#1#2#3#4#5{%
	  \reset@font\fontsize{#1}{#2pt}%
	  \fontfamily{#3}\fontseries{#4}\fontshape{#5}%
	  \selectfont}%
	  \fi\endgroup%
	  \begin{picture}(6597,2949)(4489,-1873)
		\put(11071,-556){\makebox(0,0)[lb]{\smash{{\SetFigFont{10}{12.0}{\rmdefault}{\mddefault}{\updefault}{$V$}%
		}}}}
		\put(11071,-1366){\makebox(0,0)[lb]{\smash{{\SetFigFont{10}{12.0}{\rmdefault}{\mddefault}{\updefault}{$W$}%
		}}}}
		\put(10126,-1366){\makebox(0,0)[lb]{\smash{{\SetFigFont{10}{12.0}{\rmdefault}{\mddefault}{\updefault}{$-W'$}%
		}}}}
		\put(7696,-1321){\makebox(0,0)[lb]{\smash{{\SetFigFont{10}{12.0}{\rmdefault}{\mddefault}{\updefault}{$W_1$}%
		}}}}
		\put(6751,-1321){\makebox(0,0)[lb]{\smash{{\SetFigFont{10}{12.0}{\rmdefault}{\mddefault}{\updefault}{${-}W_0$}%
		}}}}
		\put(6661,119){\makebox(0,0)[lb]{\smash{{\SetFigFont{10}{12.0}{\rmdefault}{\mddefault}{\updefault}{$I
		\times V$}%
		}}}}
		\put(7696,-556){\makebox(0,0)[lb]{\smash{{\SetFigFont{10}{12.0}{\rmdefault}{\mddefault}{\updefault}{$V_1$}%
		}}}}
		\put(7201,929){\makebox(0,0)[lb]{\smash{{\SetFigFont{10}{12.0}{\rmdefault}{\mddefault}{\updefault}{$I
		\times W$}%
		}}}}
		\put(10711,929){\makebox(0,0)[lb]{\smash{{\SetFigFont{10}{12.0}{\rmdefault}{\mddefault}{\updefault}{$Z$}%
		}}}}
		\put(6526,-421){\makebox(0,0)[lb]{\smash{{\SetFigFont{10}{12.0}{\rmdefault}{\mddefault}{\updefault}{$\text{gluing}$}%
		}}}}
		\put(5266,-1141){\makebox(0,0)[lb]{\smash{{\SetFigFont{10}{12.0}{\rmdefault}{\mddefault}{\updefault}{$D(\nu_V)$}%
		}}}}
		\put(5266,-151){\makebox(0,0)[lb]{\smash{{\SetFigFont{10}{12.0}{\rmdefault}{\mddefault}{\updefault}{$D(E)$}%
		}}}}
		\put(9688,-485){\makebox(0,0)[lb]{\smash{{\SetFigFont{10}{12.0}{\rmdefault}{\mddefault}{\updefault}{$X$}%
		}}}}
		\put(5319,-485){\makebox(0,0)[lb]{\smash{{\SetFigFont{10}{12.0}{\rmdefault}{\mddefault}{\updefault}{$X$}%
		}}}}
	  \end{picture}%
	  %%%%%%%%%%%%%%%%%%%%%%%%%%%%%%%%%%%%%%%%%%%%%%%%%%%%%%%%%%%%%%%%%%%%%%%
	  \caption{Gluing $D(E)$ and $I \x W$, and the obtained manifold pair $(Z,X)$}
	  \label{fig:gluing}
	\end{center}
  \end{figure}

  The $4$--submanifold
  \begin{equation*}
	X \cup_{V_{0}} \left( I \x V \right) \subset Z
  \end{equation*}
  (where $X$ is identified with the image of the zero--section of $E$ so that $X \subset D(E)$)
  is properly embedded in $Z$, and is bounded by $V_{1}$.
  We will rewrite $X \cup_{V_{0}} \left( I \x V \right)$ as $X$ and identify $\bd{Z}$ with $W \amalg (-W')$,
  so that
  \begin{equation*}
	\bd{(Z,X)} = (W,V) \amalg (-W',\es).
  \end{equation*}

  Now, all that is left to do is to show the existence of an $e$-class of
  $(Z,X)$ restricting to $e$ and $e'$ on the boundary components.
  Since the inclusion $W' \hookrightarrow Z \setminus X$ is homotopy
  equivalence, there exists a cohomology class $\Te \in H^{2}(Z \setminus X;\Q)$ of
  $(Z,X)$ such that $\Te|_{W'} = e'$.
  By construction, $\Te$ is an $e$-class of $(Z,X)$ and $\Te|_{W \setminus V} = e$.
  Hence, we obtain a $7$--dimensional $e$-manifold $ \b = (Z,X,\Te)$ bounded by
  \begin{equation*}
	\bd{\b} = (W,V,\Te|_{W \setminus V}) \amalg (-W',\es,\Te|_{W'})
	= \a \amalg (-\a').
    \proved
  \end{equation*}
\end{proof}

%-------------------------------------------------------------------------------------------
\subsection{\texorpdfstring{Proof of \fullref{thm:cob}}{Proof of Theorem \ref{thm:cob}}}
\label{sec:Iso-O}
%-------------------------------------------------------------------------------------------
We define a homomorphism
\begin{equation*}
  \label{eq:Phi}
  \Phi \co \O_{6}^{e} \to (\Q/\Z)^{\oplus 2}
\end{equation*}
as follows.
By \fullref{prop:pi-onto}, any element in $\O_{6}^{e}$ is represented by an
$e$-manifold of the form $(W,\es,e)$, where $W$ is a closed $6$--manifold and
$e \in H^{2}(W;\Q)$.  We then define
\begin{equation*}
  \begin{split}
	\Phi([W,\es,e]) &= \chi' ([W,e]) \\
	&\equiv \left( \frac{1}{6}\int_{W}^{}p_{1}(TW)\,e - e^{3},\
	\frac{1}{2}\int_{W}^{}e^{3} \right)   \mod{\Z^{\oplus 2}}.
  \end{split}
\end{equation*}
The rest of this section is devoted to proving that $\Phi$ is an isomorphism.
The first thing we have to do is to show that $\Phi([W,\es,e])$ is independent of the representative
$(W,\es,e)$ of $[W,\es,e]$.
\begin{lemma}
  \label{lem:Phi}
  The homomorphism $\Phi \co \O_{6}^{e} \to (\Q/\Z)^{\oplus 2}$ is well--defined.
\end{lemma}
\begin{proof}
  Let $\a= (W,\es,e)$ and $\a' = (W',\es,e')$ be any $6$--dimensional closed $e$-manifolds
  representing the same cobordism class in $\O_{6}^{e}$, and
  we prove that the difference
  $\chi([W,\es,e]) - \chi([W',\es,e'])$ belongs to $\Z^{\oplus 2}$, which implies $\chi'([W,\es,e]) =
  \chi'([W',\es,e'])$.

  There exists a $7$--dimensional $e$-manifold $\b = (Z,X,\Te)$ such that $\bd{\b} \cong \a
  \amalg (-\a')$, in particular,  
  $X$ is closed and embedded in the interior of $Z$.
  Thus, the manifold $Z_{X}$ (see \eqref{eq:ZX}) has the smooth boundary
  \begin{equation*}
	\bd{Z_{X}} = \bd{Z} \amalg (-\HX) \cong W \amalg (-W') \amalg (-\HX).
  \end{equation*}
  Since $\Te|_{\HX} = e(F_{X})$, we can write
  \begin{equation*}
	\bd{(Z_{X},\Te|_{Z_{X}})} \cong
	(W,e) \amalg (-W',e') \amalg (-\HX,e(F_{X})),
  \end{equation*}
  and this implies $[W,e] - [W',e'] = [\HX,e(F_{X})]$ in $\O_{6}(K(\Q,2))$.  
  We have
  \begin{equation*}
	\chi([W,e]) - \chi([W',e']) = \left( \Sign X, \int_{X}^{}p_{1}(\nu_{X}) \right)
  \end{equation*}
  by \fullref{lem:cu},
  and the right--hand side belongs to $\Z^{\oplus 2}$.
\end{proof}
The following theorem is the goal of this section.
\begin{theorem}
  \label{thm:O}
  The homomorphism $\Phi \co \O_{6}^{e} \to (\Q/\Z)^{\oplus 2}$ is an
  isomorphism.
\end{theorem}
\begin{proof}
  \label{proof:thmO}
  Consider the following commutative diagram:
  \begin{equation*}
	\begin{CD}
	  0 @>>>
	  \O_{4}(BSO(3))
	  @>{\up}>>
	  \O_{6}(K(\Q,2))
	  @>{\pi}>>
	  \O_{6}^{e}
	  @>>> 0\\
	  @.   @| @| @V{\Phi}VV @. \\
	  0 @>>>
	  \O_{4}(BSO(3))
	  @>{\up}>>
	  \O_{6}(K(\Q,2))
	  @>{\chi'}>>
	  (\Q/\Z)^{\oplus 2}
	  @>>> 0\\
	\end{CD}
  \end{equation*}
  The lower horizontal sequence is exact by \fullref{prop:exact-O}, and
  the homomorphism $\pi$ is surjective by \fullref{prop:pi-onto}.
  To complete the proof, we only have to show that the upper horizontal sequence is exact,
  more specifically, $\Image \up = \Ker \pi$.
  We prove this in two steps as follows.

  \textit{\hypertarget{claim:sub}{Claim\,1}\textup{:} $\Image \up \subset \Ker \pi$.}\
  Let $[X,E] \in \O_{4}(BSO(3))$ be any element, then $\pi(\up([X,E])) = [S(E),\es,e(F_{E})]$.
  Regard $X$ as the image of the zero--section of $E$ so that
  $X \subset \Int D(E)$.
  The cohomology class $e(F_{E})$ is an $e$-class of $(S(E),\es) = \bd{(D(E),X)}$, and it 
  uniquely extends to an $e$-class $e_{E}$ of $(D(E),X)$.  The obtained
  $e$-manifold $(D(E),X,e_{E})$ %\footnote{This $e$-manifold is also used in
  %\fullref{ex:es} to obtain an explicit formula of $\s$ in some special cases.}
  is bounded by $(S(E),\es,e(F_{E}))$, and hence, we have $\pi(\up([X,E])) = 0$.

  \textit{\hypertarget{claim:sup}{Claim\,2}\textup{:} $\Image \up \supset \Ker \pi$.}
  Next, we prove the opposite inclusion.
  Let $[W,e] \in \Ker \pi$ be any element, then $\a = (W,\es,e)$ bounds a
  $7$--dimensional $e$-manifold $\b = (Z,X,\Te)$, namely $\bd{\b} \cong \a$.
  In particular, we have $\Te|_{\HX} = e(F_{X})$.
  Since \[
	\bd{(Z_{X}, \Te|_{Z_{X}})} \cong (W,e) \amalg (-\HX, e(F_{X})),
  \]
  we have
  \begin{equation*}
	[W,e] = [\HX,e(F_{X})] = \up([X,\nu_{X}])
  \end{equation*}
  in $\O_{6}(K(\Q,2))$, where $\nu_{X}$ is the normal bundle of $X$.
  Therefore, $[W,e]$ belongs to $\Image \up$.
  This completes the proof.
\end{proof}
\fullref{thm:cob} can be proved very easily from \fullref{thm:O} as follows.
\begin{proof}[\hypertarget{proof:cob}{Proof of \fullref{thm:cob}}]
  For a $6$--dimensional closed $e$-manifold $\a$, its cobordism class $[\a]  \in \O_{6}^{e}
  \cong (\Q/\Z)^{\oplus 2}$ (\fullref{thm:O}) has a finite order, say $m$.  The meaning of $m[\a] = 0$  is that
  there exists a $7$--dimensional $e$-manifold $\b$ such that $\bd{\b} \cong 
  \amalg^{m}\a$.
\end{proof}

%%%%%%%%%%%%%%%%%%%%%%%%%%%%%%%%%%%%%%%%%%%%%%%%%%%%%%%%%%%%%%%%%%%%%%%%%%%%%%%%%%%%%%%%%%%%%%
\section{\texorpdfstring{Proof of \fullref{thm:main1}}{Proof of Theorem \ref{thm:main1}}}
\label{sec:proof1}
%%%%%%%%%%%%%%%%%%%%%%%%%%%%%%%%%%%%%%%%%%%%%%%%%%%%%%%%%%%%%%%%%%%%%%%%%%%%%%%%%%%%%%%%%%%%%%
In this section, we give a proof of \fullref{thm:main1}.
%-------------------------------------------------------------------------------------------
\subsection{\texorpdfstring{Signature of $4$--manifolds in $7$--manifolds}{Signature of
4-manifolds in 7-manifolds}}
\label{sec:7e}
%-------------------------------------------------------------------------------------------

For a pair $(Z,X)$, let us consider the Mayer--Vietoris exact sequence
of $(Z;Z \setminus X, U_{X})$:
\begin{multline}
  \label{eq:MV}
  \dotsb \to
  H^{2}(Z;\Q) \xrightarrow{(j_{X}^*, i_{X}^*)}
  H^{2}(Z \setminus X;\Q) \oplus H^{2}(X;\Q) \\
  \xrightarrow{\i_{X}^* - \rho_{X}^*}
  H^{2}(\HX;\Q)
  \xrightarrow{\d^*}
  H^{3}(Z;\Q)
  \to \dotsb
\end{multline}
Here, we identified $H^{2}(U_{X};\Q)$ with $H^{2}(X;\Q)$, and
$H^{2}(U_{X} \setminus X;\Q)$ with $H^{2}(\HX;\Q)$ via the isomorphisms given by the
homotopy equivalences $U_{X} \simeq X$ and $U_{X} \setminus X \simeq \HX$,
and here, the maps
\begin{equation*}
  i_{X} \co X \hookrightarrow Z,\qquad
  \i_{X} \co \HX \hookrightarrow Z \setminus X, \qquad
  j_{X} \co Z \setminus X \hookrightarrow Z
\end{equation*}
denote the inclusions.
Denote by 
\[
	t_{X} \in H^{3}(Z;\Q)
\]
the fundamental cohomology class
of $X$ (the Poincar\'e
dual of the fundamental homology class $[\pair{X}] \in H_{\dim X}^{}(\pair{Z}; \Q)$).
Since $\rho_{X}{}_{!}e(F_{X}) = 2$ (see \eqref{eq:gysin}), we have
\begin{equation}
  \label{eq:eX}
  \d^*e(F_{X}) = 2 t_{X}.
\end{equation}
The following lemma states that the existence of a quasi $e$-class of $(Z,X)$
is almost equivalent to $[\pair{X}]=0$ (exactly equivalent if the second betti--number of $\bd{X}$ vanishes).
\begin{lemma}
  \label{lem:X0-exist}
  Let $(Z,X)$ be a pair of manifolds of codimension $3$.
  Then,  $[\pair{X}] = 0$ if, and only if,
  there exist cohomology classes $e \in H^{2}(Z \setminus X;\Q)$ and $\g \in H^{2}(X;\Q)$ such that
  \begin{equation*}
	e|_{\HX} = e(F_{X}) + 2 \rho_{X}^* \g.
  \end{equation*}
  Moreover, if such $e$ and $\g$ exist and $\g|_{\bd{X}} = 0$, then $e$ is a quasi
  $e$-class of $(Z,X)$,
  and $\g$ is the self--linking form of the quasi $e$-manifold $(Z,X,e)$.
\end{lemma}
\begin{proof}
  The vanishing of $[\pair{X}]$ implies $e(F_{X}) \in \Ker \d^*$ by \eqref{eq:eX}, and thus, 
  $
  e(F_{X}) = e|_{\HX} - 2 \rho_{X}^* \g
  $
  holds
  for some elements $e \in H^{2}(Z \setminus X;\Q)$ and $\g \in H^{2}(X;\Q)$ by
  the exactness of the sequence \eqref{eq:MV}.
  The converse also holds.
  The second half of the statement is obvious from \fullref{def:eclass}\,\eqref{item:qe}.
\end{proof}
\begin{lemma}
  \label{lem:e3}
  Let $(Z,X)$ be a pair of closed manifolds of dimensions $7$ and $4$.
  If $X$ is rationally null--homologous in $Z$, then
  \begin{equation*}
	\chi([\HX, e(F_{X})]) = (\Sign X, -3 \Sign X).
  \end{equation*}
\end{lemma}
\begin{proof}
  By \fullref{lem:cu}, we have
  \[
  \chi([\HX,e(F_{X})]) = \left( \Sign X, \int_{X}^{}p_{1}(\nu_{X}) \right).
  \]
  The Hirzebruch signature theorem \eqref{eq:hirzebruch} and the vanishing of the homology class of $X$ imply
  \begin{equation*}
	\int_{X}^{}p_{1}(\nu_{X}) =
	\int_{X}^{}p_{1}(TZ) - \int_{X}^{} p_{1}(TX)
	= -\int_{X}^{}p_{1}(TX) = - 3 \Sign X.
	\proved
  \end{equation*}
\end{proof}
\begin{proposition}
  \label{prop:SL}
  Let $(Z,X)$ be a pair of closed manifolds of dimensions $7$ and $4$.  If $e$ is a quasi $e$-class of
  $(Z,X)$, then we have
  \[
  \Sign X = 4\, \L(Z,X,e).
  \]
\end{proposition}
\begin{proof}
  Let $\g \in H^{2}(X;\Q)$ be  the self--linking form of $\b = (Z,X,e)$, namely $e|_{\HX} =
  e(F_{X}) + 2\, \rho_{X}^*\g$. 
  Then, we can write
  \begin{equation*}
	\begin{split}
	  e^{3}|_{\HX} &= e(F_{X})^{3}
	  + 6 e(F_{X})^{2} \rho_{X}^*\g
	  + 12 e(F_{X})\rho_{X}^* \g^{2} + 8 \rho_{X}^*\g^{3}\\
	  &=
	  e(F_{X})^{3}
	  + 6 \rho_{X}^* \left( p_{1}(\nu_{X}) \g \right)
	  + 12 e(F_{X}) \rho_{X}^* \g^{2}
	  + 8 \rho_{X}^*\g^{3}.
	\end{split}
  \end{equation*}
  Here, we used the relation \eqref{eq:pe1} (where $E = \nu_{X}$)
  on the second term of the right--hand side.
  Since $\dim X < 6$, 
  the second and the last terms of the right--hand side vanish.
  Integrating the both sides over $\HX$, we obtain
  \begin{align*}
	\int_{\HX}^{} e^{3}
	%		&= \int_{\HX}^{}e(F_{X})^{3} + 12 e(F_{X}) \rho_{X}^*\g^{2}\\
	&= 2 \chi_{2}(\HX,e(F_{X})) + 24 \L(\b).
  \end{align*}
  The left--hand side vanishes by Stokes' theorem, since $\bd{(Z_{X}, e^{3}|_{Z_{X}})} =
  (-\HX,e^{3}|_{\HX})$.
  Thus, we obtain $\chi_{2}(\HX,e(F_{X})) = - 12 \L(\b)$.

  On the other hand, the existence of a quasi $e$-class of $(Z,X)$ implies that $X$ is rationally
  null--homologous in $Z$ by \fullref{lem:X0-exist}, and so $(Z,X)$ satisfies the assumption of
  \fullref{lem:e3}.
  Hence, we obtain
  \begin{equation*}
	\Sign X = -\frac{1}{3}\chi_{2}(\HX, e(F_{X})) = 4 \L(\b).
	\proved
  \end{equation*}
\end{proof}
As a corollary, we obtain the following vanishing property of the signature.
\begin{corollary}
  \label{cor:O4-X}
  If a pair $(Z,X)$ of closed manifolds of dimensions $7$ and $4$ admits an $e$-class, then
  $\Sign X = 0$.
\end{corollary}
\begin{proof}
  Let $e$ be an $e$-class of $(Z,X)$, then $\L(Z,X,e) = 0$ by the definition of
  $\L$.  This implies $\Sign X = 0$ by \fullref{prop:SL}.
\end{proof}

%-------------------------------------------------------------------------------------------
\subsection{\texorpdfstring{Definition of the invariant $\s(\a)$}{Definition of the invariant s}}
\label{sec:def-s}
%-------------------------------------------------------------------------------------------
Let $\a$ be a $6$--dimensional closed $e$-manifold.
We first review and generalize the definition of the invariant $\s(\a)$ as follows.
By \fullref{thm:cob}, there exists a $7$--dimensional $e$-manifold $\b =
(Z,X,\Te)$ such that $ \bd{\b} \cong \amalg^{m} \a$ for some positive integer $m$.
More generally, we shall assume that $\b$ is a quasi $e$-manifold.
We then define
\begin{equation}
  \label{eq:def-s}
  \s(\a) = \frac{\Sign X - 4 \L(\b)}{m} \in \Q.
\end{equation}
In this section, we use this definition instead of the one given in \fullref{thm:main1}, because 
\eqref{eq:def-s} includes the formula in \fullref{thm:main1}\,\eqref{item:sl}.
\begin{proposition}\label{prop:s}
  For a $6$--dimensional closed $e$-manifold $\a$,
  the rational number $\s(\a)$, defined as in \eqref{eq:def-s}, depends only on the isomorphism class of $\a$.
\end{proposition}
\begin{proof}
  Let $\calN$ be the set of the isomorphism classes of $6$--dimensional null--cobordant closed
  $e$-manifolds.
  Note that if $\a \in \calN$, then $\amalg^{m}\a \in \calN$ for any positive
  integer $m$.
  We prove the statement in three steps as follows.
  The first step is the most important, and the rests are proved in a formal way.

  \textit{\hypertarget{claim:N}{Claim\,1}\textup{:} If $\a \in \calN$, then $\s(\a)$ is well--defined.}
  Assume $\a \in \calN$, and take any $7$--dimensional quasi $e$-manifolds $\b_{0} =
  (Z_{0},X_{0},e_{0})$ and $\b_{1} =
  (Z_{1},X_{1}, e_{1})$
  equipped with fixed identifications $\bd{\b_{0}} = \a = \bd{\b_{1}}$. We need to show the equality
  \begin{equation}
	\label{eq:SL}
	\Sign X_{0} -4 \L(\b_{0}) = \Sign X_{1} - 4 \L(\b_{1}).
  \end{equation}
  For that, we consider the pair of closed manifolds
  \begin{equation*}
	(Z,X) = \left(
	Z_{0} \cup (-Z_{1}),
	X_{0} \cup (-X_{1})
	\right)
  \end{equation*}
  obtained from $(Z_{0},X_{0})$ and $(-Z_{1},-X_{1})$ by gluing along the identified
  boundaries, and write $(W,V) = (\bd{Z_{0},\bd{X_{0}}}) \subset (Z,X)$.
  Since the quasi $e$-classes $e_{0}$ and $e_{1}$ restricts to the same $e$-class of $(W,V)$,
  there exists a quasi $e$-class $e$ of $(Z,X)$ such
  that $e|_{Z_{i} \setminus X_{i}} = e_{i}$ for $i=0,1$.
  Thus, we obtain a $7$--dimensional closed quasi $e$-manifold $\b = (Z,X,e)$.

  The self--linking form $\g \in H^{2}(X;\Q)$ of $\b$ is trivial on $V$, and thus,
  there exists an element
  \[
  \Tg = (\Tg_{0}, \Tg_{1}) \in H^{2}(X,V;\Q) =
  H^{2}(\pair{X_{0}};\Q)
  \oplus
  H^{2}(\pair{X_{1}};\Q)
  \]
  such that the homomorphism $H^{2}(X,V;\Q) \to H^{2}(X;\Q)$ maps $\Tg$ to $\g$.
  Similarly, the homomorphism $H^{2}(\pair{X_{i}};\Q) \to H^{2}(X_{i};\Q)$ maps
  $\Tg_{i}$ to the self--linking form of $\b_{i}$, and thus
  \begin{equation*}
	\L(\b_{i}) = \int_{X_{i}}^{}\Tg_{i}^{2}.
  \end{equation*}
  Hence, we have
  \begin{equation*}
	\L(\b) = \int_{X}^{} \g^{2}
	=
	\int_{X_{0}}^{}\Tg_{0}^{2} - \int_{X_{1}}^{}\Tg_{1}^{2}
	= \L(\b_{0}) - \L(\b_{1}).
  \end{equation*}
  (Namely, the self--linking number $\L$ is additive with respect to the decompositions of closed quasi
  $e$-manifolds.)
  By the additive properties of $\L$ and the signature, and by  \fullref{prop:SL}, we have
  $\Sign X_{0} - \Sign X_{1} = 4 \L(\b_{0}) - 4 \L(\b_{1})$, and so \eqref{eq:SL} holds.

  \textit{\hypertarget{claim:additive}{Claim\,2}\textup{:} If $\a \in \calN$, then
  $\s(\amalg^{m} \a) = m \s(\a)$ for any positive integer $m$.}
  Let $\b = (Z,X,e)$ be a $7$--dimensional $e$-manifold such that $\bd{\b} \cong \a$, then
  $\amalg^{m} \a$ bounds $\amalg^{m}\b$.
  The rational number $\s(\amalg^{m} \a)$ is well--defined by \hyperlink{claim:N}{Claim\,1}, and we have
  \begin{equation*}
	\s(\amalg^{m} \a) = \Sign (\amalg^{m} X) = m\,\Sign X = m\,\s(\a).
  \end{equation*}

  \textit{\hypertarget{claim:all}{Claim\,3}\textup{:} $\s(\a)$ is well--defined for any
  $6$--dimensional closed $e$-manifold $\a$.}
  Let $\a$ be a closed $e$-manifold, and
  $m$ a positive integer such that $\amalg^{m} \a \in \calN$
  (such $m$ exists by \fullref{thm:cob}).
  The rational number $\s(\amalg^{m} \a)$ is well--defined by \hyperlink{claim:N}{Claim\,1}.
  We can show that the rational number $\s(\amalg^{m}\a)/m$ does not depend on the choice of $m$ as follows.

  If $\amalg^{m} \a \in \calN$ and $\amalg^{m'} \a \in \calN$ for some positive integers $m$ and $m'$, then 
  $\amalg^{m m'} \a \in \calN$.
  Thus, the rational numbers $\s(\amalg^{m} \a)$, $\s(\amalg^{m'} \a)$, and
  $\s(\amalg^{m m'}\a)$ are well--defined by \hyperlink{claim:N}{Claim\,1}.
  Since $\amalg^{m'}(\amalg^{m}\a) = \amalg^{m m'}\a = \amalg^{m}(\amalg^{m'} \a)$, we have
  \begin{equation*}
	\frac{\s(\amalg^{m}\a)}{m} =
	\frac{\s(\amalg^{m m'}\a)}{m m'} =
	\frac{\s(\amalg^{m'}\a)}{m'}
  \end{equation*}
  by \hyperlink{claim:additive}{Claim\,2}. This implies that $\s(\amalg^{m} \a)/m$ does not depend on
  the choice of $m$.
\end{proof}
%-------------------------------------------------------------------------------------------
\subsection{\texorpdfstring{Proof of \fullref{thm:main1}}{Proof of Theorem \ref{thm:main1}}}
\label{sec:Proof-2}
%-------------------------------------------------------------------------------------------
By using the results we have obtained so far, we prove \fullref{thm:main1}. 
\begin{proof}[Proof of \fullref{thm:main1}]
  By \fullref{prop:s}, the rational number $\s(\a)$ is well--defined for any $6$--dimensional
  closed $e$-manifold $\a$.
  \fullref{axiom2} and \fullref{thm:main1}\,\eqref{item:sl} are obvious from the definition
  \eqref{eq:def-s}.

  We prove that $\s$ satisfies \fullref{axiom1} as follows.
  Let $\a$ be a $6$--dimensional closed $e$-manifold, and $\b = (Z,X,e)$ a $7$--dimensional
  $e$-manifold such that $\bd{\b} \cong \amalg^{m} \a$ for some positive integer $m$.
  Then, $\bd{(-\b)} \cong \amalg^{m}(-\a)$, and the definition of $\s$ implies
  \begin{equation*}
	\s(-\a) = \frac{\Sign (-X)}{m} =  -\s(\a).
  \end{equation*}
  Next, let $\a'$ be another $6$--dimensional closed $e$-manifold, and $\b' = (Z',X',e')$ a
  $7$--dimensional $e$-manifold such that $\bd{\b'} \cong \amalg^{m'}\a$ for some positive
  integer $m'$.
  Then $(\amalg^{m'}\b) \amalg (\amalg^{m}\b')$ bounds
  $\amalg^{m m'} (\a \amalg \a')$, and we have
  \begin{align*} \s(\a \amalg \a') &= \frac{m' \Sign X + m
	\Sign X'}{m m'}\\ &= \frac{\Sign X}{m} + \frac{\Sign X'}{m'}\\ &= \s(\a) + \s(\a').
  \end{align*}
  Hence, \fullref{axiom1} holds.

  We prove \fullref{thm:main1}\,\eqref{item:unique} (uniqueness of $\s$) as follows.
  Let $\s'$ be an invariant of the isomorphism classes of
  $6$--dimensional closed $e$-manifolds satisfying the axioms.
  Let us consider the difference
  \[
  f(\a) = \s'(\a) - \s(\a) \in \Q.
  \]
  If two $6$--dimensional closed $e$-manifolds $\a_{0}$ and $\a_{1}$ are cobordant, that is, if
  there exists a $7$--dimensional 
  $e$-manifold $\b = (Z,X,e)$  such that $\bd{\b} \cong \a_{0} \amalg (-\a_{1})$, then
  \begin{equation*}
	f(\a_{0}) - f(\a_{1}) = \s'(\bd{\b}) - \s(\bd{\b}) = \Sign X - \Sign X = 0
  \end{equation*}
  by the axioms.  
  Thus, we can regard $f$ as a function on $\O_{6}^{e}$:
  \begin{equation*}
	f \co \O_{6}^{e} \to \Q
  \end{equation*}
  Moreover, \fullref{axiom1} implies that $f$ is a homomorphism.

  On the other hand,
  any homomorphism $\O_{6}^{e} \to \Q$ is trivial by \fullref{thm:cob}, in particular,
  $f$ must be trivial.   Namely, $\s' = \s$.
\end{proof}
This completes the proof of \fullref{thm:main1}.

%%%%%%%%%%%%%%%%%%%%%%%%%%%%%%%%%%%%%%%%%%%%%%%%%%%%%%%%%%%%%%%%%%%%%%%%%%%%%%%%%%%%%%%%%%%%%
\section{\texorpdfstring{Existence and uniqueness of $e$-classes}{Existence and uniqueness of
e-classes}}
\label{sec:e}
%%%%%%%%%%%%%%%%%%%%%%%%%%%%%%%%%%%%%%%%%%%%%%%%%%%%%%%%%%%%%%%%%%%%%%%%%%%%%%%%%%%%%%%%%%%%%
As in \fullref{sec:emb},
we denote by $\calE_{Z,X} \subset H^{2}(Z \setminus X; \Q)$ the set of all $e$-classes of a
manifold pair $(Z,X)$ of codimension $3$. 
In this short section, we study some elementary properties of $\calE_{Z,X}$, and then we give
a proof of \fullref{thm:Q}.

For elements $e \in \calE_{Z,X}$  and $a \in \Ker i_{X}^{*}$, 
where 
\[
	i_{X}^* \co H^{2}(Z;\Q) \to H^{2}(X;\Q)
\]
is the restriction, the cohomology class 
\[
g_{e}(a) \defeq e + a|_{Z \setminus X} \in H_{}^{2}(Z \setminus X; \Q)
\]
belongs to $\calE_{Z,X}$, because $a|_{\HX} = 0$.
Thus, we obtain an affine homomorphism
\begin{equation*}
  g_{e} \co \Ker i_{X}^{*} \to \calE_{Z,X}.
\end{equation*}
Note that $g_{e}$ is defined only when $\calE_{Z,X} \neq \es$.
Let $\#A \in \N \cup \left\{ \infty \right\}$ denote the number of elements in
a set $A$, and let
$[\pair{X}] \in H_{\dim X}(\pair{Z};\Q)$ be the fundamental homology class of
$(\pair{X})$.
\begin{proposition}\label{prop:eclass}
  The following statements hold\textup{:}
  \begin{enumerate}
	\item \label{item:quasi} If $[\pair{X}] \neq 0$, then $\calE_{Z,X} = \es$.
	\item\label{item:bij} For any $e \in \calE_{Z,X}$, the map
	  $g_{e} \co \Ker i_{X}^{*} \to \calE_{Z,X}$ is an affine isomorphism.
	\item\label{item:onto} Assume that $i_{X}^*$ is surjective.  Then,
	  $\calE_{Z,X} \neq \es$ if, and only if, $[\pair{X}] = 0$.
	\item\label{item:inj}
	  $\#\calE_{Z,X} = 1$ if, and only if, $(Z,X)$ is simple.
	\item\label{item:iso}
	  \textup{\fullref{prop:simple}} holds.  Namely, when $i_{X}^*$ is an isomorphism,
	  $(Z,X)$ is simple if, and only if, $[\pair{X}] = 0$.
  \end{enumerate}
\end{proposition}
\begin{proof}
  \eqref{item:quasi} is a direct consequence of \fullref{lem:X0-exist}.  In fact, if $[\pair{X}] \neq
  0$, then $(Z,X)$ does not even admit a quasi $e$-class.

  \eqref{item:bij}\
  Since  two elements in $\calE_{Z,X}$ differ by an element in the kernel $\Ker \i_{X}^*$ of
  the homomorphism $\i_{X}^{*} \co H^{2}(Z \setminus X;\Q) \to H^{2}(\HX;\Q)$, we only need
  to show that the homomorphism
  $
	j_{X}^* \co H^{2}(Z;\Q) \to H^{2}(Z \setminus X;\Q)
  $
  restricts a bijection $\Ker i_{X}^* \to \Ker \i_{X}^*$, and this is immediate from the
  following commutative diagram: 
  \begin{equation*}
	\xymatrix{
	H^{2}(Z\setminus X, \HX;\Q)  \ar[r]
	&H^{2}(Z \setminus X;\Q)  \ar[r]^{\i_{X}^{*}}
	&H^{2}(\HX;\Q)\\
	H^{2}(Z,X;\Q)  \ar[u]^{f}_{\cong} \ar[r]
	&H^{2}(Z;\Q)  \ar@{^{(}->}[u]^{j_{X}^*} \ar[r]^{i_{X}^{*}}
	&H^{2}(X;\Q)
	}
  \end{equation*}
  Here, the horizontal sequences are the long exact sequences of $(Z \setminus X,
  \HX)$ and $(Z,X)$, and $f$ is the excision isomorphism.
  Note that the vanishing of $H^{2}(Z, Z \setminus X;\Q)$ follows the injectivity of
  $j_{X}^{*}$.

  \eqref{item:onto}\ We only give a proof of that the vanishing of $[\pair{X}]$ implies
  $\calE_{Z,X} \neq \es$, since
  the converse is given by \eqref{item:quasi}.
  Assume that $i_{X}^*$ is surjective and $[\pair{X}] = 0$.  By \fullref{lem:X0-exist}, there exist
  $e \in H^{2}(Z\setminus X;\Q)$ and $\g \in H^{2}(X;\Q)$ such that
  $e|_{\HX} = e(F_{X}) + 2 \rho_{X}^*\g$.
  By the assumption, there exists $\varepsilon \in H^{2}(Z;\Q)$ such that
  $\varepsilon|_{X}=\g$. The homotopy equivalence $j_{X} \i_{X} \simeq i_{X} \rho_{X} \co \HX
  \to Z$ implies
  $\varepsilon|_{\HX} = \rho_{X}^*\g$, and
  the cohomology class $e' = e - 2 \varepsilon|_{Z \setminus X} \in H^{2}(Z \setminus X;\Q)$ satisfies
  \begin{equation*}
	e'|_{\HX} =
	\left( e(F_{X}) + 2 \rho_{X}^*\g \right)
	-
	2 \rho_{X}^*\g
	=
	e(F_{X}).
  \end{equation*}
  This means $e' \in \calE_{Z,X} \neq \es$.

  \eqref{item:inj}\ By \eqref{item:bij}, $\#\calE_{Z,X} = 1$ implies $\Ker i_{X}^{*} = \left\{ 0
  \right\}$ which means that $i_{X}^*$ is injective, and thus, $(Z,X)$ is simple.
  Conversely, let us assume that $(Z,X)$ is simple, that is, $i_{X}^*$ is injective and $\calE_{Z,X} \neq \es$.
  Fix any element $e \in \calE_{Z,X}$, then the map $g_{e} \co
  \Ker i_{X}^{*} = \left\{ 0 \right\} \to \calE_{Z,X}$ is a bijection by \eqref{item:bij}, and thus,
  $\#\calE_{Z,X} = \#\left\{ e \right\} = 1$.

  \eqref{item:iso}\ This is a combination of \eqref{item:onto} and \eqref{item:inj}.
\end{proof}
If we drop the surjectivity assumption of $i_{X}^*$ in \fullref{prop:eclass}\,\eqref{item:onto}, then the
statement does not hold anymore, and a counterexample is given in the following.
\begin{remark}
  \label{rem:S7}
  It is known that any oriented closed smooth $4$--manifold $X$ can be smoothly embedded in $S^{7}$
  (c.f.~\cite[Theorem\,9.1.23, Remark\,9.1.24]{gompf-stipsicz}).
  Let us assume $\Sign X \neq 0$ (for example, $\Sign \CP^{2} = 1$).
  Obviously, $X$ is null--homologous in $S^{7}$, but $(S^{7}, X)$ does not admit any
  $e$-class by \fullref{cor:O4-X}.
\end{remark}

Now, \fullref{thm:Q} is proved as follows.
\begin{proof}[\hypertarget{proof:Q}{Proof of \fullref{thm:Q}}]
  Let us assume that $(W,V)$ is simple.
  By \fullref{cor:ftn}, the image $\Image \s_{W,V} \subset \Q$ of the 
  function $\s_{W,V} \co \calE_{W,V} \to \Q$ depends only on the 
  isomorphism class of $(W,V)$.
  The simplicity of $(W,V)$ implies that $\calE_{W,V}$ consists of just one element, say $\calE_{W,V} = \left\{ e
  \right\}$, and therefore, $\Image \s_{W,V} = \left\{ \s(W,V,e) \right\}$.
  Thus, $\s(W,V,e)$ is an invariant of the isomorphism class of $(W,V)$.
\end{proof}

%%%%%%%%%%%%%%%%%%%%%%%%%%%%%%%%%%%%%%%%%%%%%%%%%%%%%%%%%%%%%%%%%%%%%%%%%%%%%%%%%%%%%%%%%%%%
\section{Seifert surfaces}
\label{sec:seifert-intro}
%%%%%%%%%%%%%%%%%%%%%%%%%%%%%%%%%%%%%%%%%%%%%%%%%%%%%%%%%%%%%%%%%%%%%%%%%%%%%%%%%%%%%%%%%%%%
In this section, we establish the relationship between Seifert surfaces and
quasi $e$--classes.
When we consider the intersection of submanifolds,
we will always assume that the submanifolds are in general positions (by deforming them
slightly if necessary) so that the intersections becomes smooth manifolds.

Let $(X,E)$ be a pair of a manifold $X$ and an oriented  vector bundle $E$ of rank $3$ over
$X$. 
We denote by
\[
	\tau \co E \to E,\quad v \mapsto -v
\]
the involution given by the multiplication by a scalar $-1$.
There is a direct sum decomposition
\begin{equation}
  \label{eq:Hpm}
  H^{2}(S(E);\Q) = H_{+1} \oplus H_{-1},
\end{equation}
where
$H_{\pm 1}$ is the eigenspace of the involution $\t^* \co H^{2}(S(E);\Q) \to H^{2}(S(E);\Q)$ with the eigenvalue $\pm 1$.
The subspace
$H_{+1}$ is the image of the pull--back $\rho_{E}^* \co H^{2}(X;\Q) \to H^{2}(S(E);\Q)$,
and $H_{-1}$ is the subspace spanned by the Euler class class $e(F_{E})$ of $F_{E}$.
These facts are proved by using the Thom--Gysin exact sequence of $E$.

The Euler class $e(F_{E})$ of $F_{E}$ is algebraically characterized as follows.  
\begin{lemma}\label{lem:ie}
  Let $(X,E)$ be as above.
  If a cohomology class $a \in H^{2}(S(E);\Q)$ satisfies two conditions
  \begin{enumerate}
	\item \label{item:tau} $\tau^*a = -a$,
	\item \label{item:g2} $\rho_{E}{}_{!}a = 2$,
  \end{enumerate}
  then $a = e(F_{E})$ over $\Q$.
  Here, $\rho_{E}{}_{!} \co H^{2}(S(E);\Q) \to H^{0}(X;\Q)$ is the Gysin homomorphism of the
  associated sphere bundle $\rho_{E} \co S(E) \to X$.
\end{lemma}
\begin{proof}
  The property \eqref{item:tau} implies $a \in H_{-}$, and \eqref{item:g2} implies $a =
  e(F_{E})$.
\end{proof}
For a manifold pair $(Z,X)$ of codimension $3$, a \emph{Seifert surface} of $X$ is a
proper (oriented) submanifold $Y$ of $Z_{X}$ (see \eqref{eq:ZX}) such that $Y \cap \HX = s(X)$ for some
section $s \co X \to \HX$, and such that
the natural isomorphism 
$\nu_{Y}|_{s(X)} \cong F_{X}|_{s(X)}$
preserves the orientation.
Note that $Y$ may have the corner $s(\bd{X})$ (see also \eqref{eq:bdry-Y} below).
More generally, if $Y \setminus U_{s(X)}$ is immersed in $Z_{X} \setminus U_{\HX}$ for some
open neighborhoods $U_{s(X)}$ of $s(X)$ and $U_{\HX}$ of $\HX$ (and so $Y$ has no multiple
points on $U_{s(X)}$), then we say $Y$ is an \emph{immersed} Seifert surface.

Let $Y$ be such an immersed Seifert surface of $X$ in $Z$.
Define
\begin{equation}
  \label{eq:Y-bdry}
  F(s) \defeq s^*F_{X} \cong s^*(\nu_{Y}|_{s(X)}),
\end{equation}
which is an oriented vector bundle of rank $2$ over $X$, and so its Euler class
\[
	e(F(s)) \in H^{2}(X;\Z)
\]
is defined. 
Unless otherwise stated, we do not assume that $F(s)$ is trivial in this paper.
Let
\[
	t_{Y} \in H^{2}(Z \setminus X;\Q) \cong H^{2}(Z_{X};\Q)
\]
be the fundamental cohomology class of $Y$, then we have
\begin{equation}	
  \label{eq:sty}
  s^*(t_{Y}|_{\HX}) = e(F(s)),\qquad
  s^*\t^* (t_{Y}|_{\HX}) = 0.
\end{equation}

Now, let us write 
\[
	(W,V) = \bd{(Z,X)},\quad S = Y \cap W_{V}
\]
for a moment.
Then,
$S$ is an immersed Seifert surface of $V$ in $W$ with respect to the section $s|_{V} \co V \to
\HV$ (namely $S \cap \HV = s(V)$).  
Note that
\begin{equation*}
  t_{S} = t_{Y}|_{W \setminus V},
  \qquad 
  F(s|_{V}) = F(s)|_{V},
  \qquad
  \nu_{S} = \nu_{Y}|_{S}
\end{equation*}
by definition.
In line with our orientation conventions, the (oriented) boundaries of 
$Y$ and $S$ are given as follows:
\begin{equation}
  \label{eq:bdry-Y}
  \bd{Y} = (\pm S) \cup (\mp s(X)),\qquad
  \bd{S} = \pm s(V)
\end{equation}
Here, $\pm = (-1)^{\dim W}$ and $\mp = (-1)^{\dim Z}$.
\begin{remark}
Let
$Y' \subset Z$ be an (immersed) submanifold with the boundary $\bd{Y'} = (\pm (\bd{Y'} \cap W))
\cup (\mp X)$ and with the 
corner $\angle Y = \bd{X}$ such that $Y'$ intersects $W$ transversely, and we 
assume that a neighborhood of $X \subset Y$ has no multiple points.
Then, $Y= Y' \cap Z_{X}$ is an (immersed) Seifert surface of $X$ in the sense described above.
In order to avoid introducing too much notation, we will also call $Y'$ an (immersed) {Seifert surface} of $X$.
The corresponding section $s \co X \to \HX$ is defined such
that $Y' \cap \HX = s(X)$, and so $F(s) \cong \nu_{Y'}|_{X}$.
The cohomology class $t_{Y} \in H^{2}(Z \setminus X;\Q)$ is nothing but the fundamental cohomology class of
$Y' \setminus X$ (the Poincar\'e dual of the locally finite
fundamental homology class $[Y' \setminus X] \in H_{\dim Y'}^{lf}(Z \setminus X;\Q)$ of $Y' \setminus X$).
\end{remark}
The following proposition states that an immersed Seifert surface implies a quasi
$e$-class.
\begin{proposition}
  \label{prop:seifert}
  Let $(Z,X)$ be a pair of manifolds of codimension $3$, and $Y$ an immersed Seifert surface
  of $X$ with respect to a section $s \co X \to \HX$ such that $e(F(s))|_{\bd{X}} = 0$ over $\Q$.
  Then, $2 t_{Y}$ is a quasi $e$-class of $(Z,X)$, and  
  $e(F(s))/2$ is the self--linking form of the quasi $e$-manifold
  $(Z,X,2t_{Y})$.
\end{proposition}
\begin{proof}
  We write $b = t_{Y}|_{\HX} \in H^{2}(\HX;\Q)$ in this proof.
  By \fullref{lem:ie}, we have 
  \begin{equation}
	\label{eq:ap}
	b - \t^*b = e(F_{X}).
  \end{equation}
  Since $b + \t^*b$ belongs to $H_{+} = \Image \rho_{X}^*$ (see \eqref{eq:Hpm}, where $E =
  \nu_{X}$),
  there exists $c \in
  H^{2}(X;\Q)$ such that
  $b + \t^*b = \rho_{X}^*c$.
  The pull--back $s^* \co H^{2}(\HX;\Q) \to
  H^{2}(X;\Q)$ is a left--inverse of $\rho_{X}^*$, and so
  \begin{equation*}
	c = s^*(\rho_{X}^*c) = s^*b + s^*\t^*b = e(F(s))
  \end{equation*}
  by \eqref{eq:sty}, and thus
  \begin{equation}
	\label{eq:am}
	b + \t^*b = \rho_{X}^*e(F(s)).
  \end{equation}
  By \eqref{eq:ap} and \eqref{eq:am}, we have
  \begin{equation*}
	2 t_{Y}|_{\HX} = e(F_{X}) + 2 \left( \rho_{X}^*e(F(s))/2 \right).
  \end{equation*}
  Since $e(F(s))|_{\bd{X}} = 0$, $2 t_{Y}$ is a quasi $e$-class of $(Z,X)$, and  $e(F(s))/2$ is the
  self--linking form of $(Z,X,2t_{Y})$. 
\end{proof}
The following is a direct consequence of \fullref{prop:seifert}.
\begin{corollary}
  \label{cor:seifert-e}
  Let $(Z,X)$ and $Y$ be as in \textup{\fullref{prop:seifert}}.  
  If $e(F(s)) = 0$, then $2 t_{Y}$ is an $e$-class of
  $(Z,X)$.   
\end{corollary}
\begin{proof}
  By \fullref{prop:seifert}, the self--linking form of
  $(Z,X,2t_{Y})$ vanishes, in other words, $2t_{Y}$ is an $e$-class by definition.
\end{proof}
%%%%%%%%%%%%%%%%%%%%%%%%%%%%%%%%%%%%%%%%%%%%%%%%%%%%%%%%%%%%%%%%%%%%%%%%%%%%%%%%%%%%
\section{\texorpdfstring{Haefliger's invariant}{Haefliger's invariant}}
\label{sec:3knot}
%%%%%%%%%%%%%%%%%%%%%%%%%%%%%%%%%%%%%%%%%%%%%%%%%%%%%%%%%%%%%%%%%%%%%%%%%%%%%%%%%%%%
In this section, we prove \fullref{thm:Hf}.
We also prove the geometric formula for $\s$, and as a corollary, we obtain 
more general results (\fullref{cor:takase})
establishing the relationship between Takase's invariant $\O$
and our invariant $\s$.

%-------------------------------------------------------------------------------------------
\subsection{Review of Haefliger's invariant}
\label{sec:review-h}
%-------------------------------------------------------------------------------------------
We begin by reviewing Haefliger's results~\cite{Haefliger-knot}~\cite{Haefliger-diff} on the
classification of smooth $3$--knots in $S^{6}$.
He first showed that the set
$\Emb(S^{3},S^{6})$ of the isotopy classes of smooth embeddings $f \co S^{3} \to S^{6}$ is an
abelian group with the group structure given by the connected sum.
Write $M_{f} = f(S^{3})$.
He showed the existence of an oriented proper framed $4$--submanifold $X \subset
D^{7}$ such that $\bd{X} = M_{f}$ and $\Sign X = 0$.
Here, a \emph{framing} of $X$ is the homotopy class of a triple $(s_{1}, s_{2}, s_{3})$, $s_{i} \co X \to
\nu_{X}$, of linearly independent sections of the normal bundle $\nu_{X}$ of $X$.
We shall assume $s_{i}(X) \subset \HX$.
The homomorphism 
\[
H^{2}(\pair{X};\Q) \xrightarrow{\cong} H^{2}(X;\Q)
\]
is an isomorphism, and we will identify these two groups.
For a $2$--cycle $c$ of $X$,
the linking number $\lk(s_{1}(c), X ) \in \Q$ of $s_{1}(c)$ with $X$ is well--defined,  and it depends only on
the homology class $[c] \in H_{2}(X;\Z)$ of $c$.
Thus, we obtain a homomorphism 
\[
	\l \co H_{2}(X;\Z) \to \Q,\quad
	\l([c]) = \lk(s_{1}(c), X),
\]
which gives a cohomology class  $\l \in H^{2}(\pair{X};\Q)$.
He proved that the integral
\begin{equation}\label{eq:H}
  H(f)= \frac{1}{2}\int_{X}^{} \l^{2}
\end{equation}
is an integer and depends only on the isotopy class of $f$,
and that the induced map $H \co \Emb(S^{3},S^{6}) \to \Z$ is an
isomorphism of abelian group.

%-------------------------------------------------------------------------------------------
\subsection{\texorpdfstring{Invariant $\s(S^{6},M_{f})$}{Invariant s(S6,Mf)}}
\label{sec:s-S}
%-------------------------------------------------------------------------------------------
Let $M_{f}$ and $X$ be as before.
The pair $(S^{6},M_{f})$ is simple by \fullref{prop:simple},
and we can define the invariant $\s(S^{6},M_{f}) \in \Q$ by \fullref{thm:Q}.
By \fullref{lem:X0-exist},
there exist $e \in H^{2}(D^{7} \setminus X;\Q)$ and $\g \in H^{2}(X;\Q)$ such that
\begin{equation*}
  e|_{\HX} = e(F_{X}) + 2 \rho_{X}^*\g.
\end{equation*}
Since $\g|_{M_{f}} = 0 \in H^{2}(M_{f};\Q) = 0$,
$e$ is a quasi $e$-class of $(D^{7},X)$, and $\g$ is the self--linking form of the quasi
$e$-manifold $(D^{7},X,e)$.
In particular, $e_{f} = e|_{S^{6} \setminus M_{f}} \in H^{2}(S^{6} \setminus M_{f};\Q)$ is the unique
$e$-class of $(S^{6},M_{f})$, and 
\[
	\bd{(D^{7},X,e)} = (S^{6}, M_{f}, e_{f})
\]
as $e$-manifolds.
By \fullref{thm:main1}\,\eqref{item:sl},
we obtain a formula
\begin{equation}\label{eq:HL2}
  \s(S^{6}, M_{f}) = - 4 \int_{X}^{}\g^2.
\end{equation}
%-------------------------------------------------------------------------------------------
\subsection{\texorpdfstring{Proof of \fullref{thm:Hf}}{Proof of Theorem \ref{thm:Hf}}}
\label{sec:proof5}
%-------------------------------------------------------------------------------------------
The essential part of the proof of \fullref{thm:Hf} is that the cohomology classes $\l$ and
$\g$, we defined in this section, are actually the same.
\begin{proposition}\label{prop:lg}
  We have $\l = \g$ in $H^{2}(\pair{X};\Q)$.
\end{proposition}
\begin{proof}
  Fix a $2$-cycle $c$ of $X$, and let $\o \in H^{2}(\pair{X};\Q)$ be the Poincar\'e dual of the
  homology class $[c] \in H_{2}(X;\Q)$ of $c$.
  Let us write $\Tc = s_{1}(c)$ and $\TX = s_{1}(X)$ which are cycles of $X$. Note that
  the Poincar\'e dual of $[\pair{\TX}] \in H_{4}({\HX};\Q)$
  is $e(F_{X})/2$.

  Since $\Tc = \TX \cap \rho_{X}^{-1}(c)$, its Poincar\'e dual is
  \[
  \frac{1}{2}{e(F_{X})\,\rho_{X}^*\o} \in H^{4}(\pair{\HX};\Q).
  \]
  The homomorphism
  \[
  H^{2}(D^{7} \setminus X;\Q) \to H^{3}(D^{7}, D^{7} \setminus X;\Q)
  \]
  given by
  the pair $(D^{7}, D^{7} \setminus X)$ maps $e/2$ to
  the Thom class of $\nu_{X}$, and thus, we have
  \begin{equation*}
	\lk(\Tc, X) =
	\bk{[\Tc], e/2} =
	\frac{1}{4}\int_{\HX}^{}  e(F_{X})\,e\,\rho_{X}^*\o.
  \end{equation*}
  Therefore, we have 
  \begin{equation*}
	\l([c]) = \frac{1}{4}\int_{\HX}^{} e(F_{X})\,e\,\rho_{X}^*\o
	= \frac{1}{2}\int_{\HX}^{}e(F_{X})\,\rho_{X}^*\left( \g\,\o \right)
	= \int_{X}^{}\g\,\o
	= \bk{[c], \g},
  \end{equation*}
  where we used the relation $e|_{\HX} = e(F_{X}) + 2 \rho_{X}^*\g$ and the vanishing
  \[
  \int_{\HX}^{}e(F_{X})^{2}\rho_{X}^*\o = \int_{\HX}^{} \rho_{X}^*\left( p_{1}(\nu_{X}) \o
  \right) = 0
  \]
  following from \eqref{eq:pe1} and $\dim X < 6$.
  Hence, $\l([c]) = \bk{[c],\g}$ for any $2$--cycle $c$ of $X$, and thus, $\l = \g \in
  H^{2}(\pair{X};\Q)$.
\end{proof}
\fullref{thm:Hf} is now quite easy to proof.
\begin{proof}[Proof of \fullref{thm:Hf}]
  By \eqref{eq:H}, \eqref{eq:HL2}, and \fullref{prop:lg}, we have
  \begin{equation*}
	\s(S^{6},M_{f}) = -4 \int_{X}^{}\g^{2} = -4\int_{X}^{}\l^{2} = -8 H(f).
	\proved
  \end{equation*}
\end{proof}
We shall say that $\s$ is a natural generalization of Haefliger's
invariant $H$.

%-------------------------------------------------------------------------------------
\subsection{Geometric formula}
\label{sec:geom-formula}
%-------------------------------------------------------------------------------------
Let $(W,V)$ be a pair of closed manifolds of dimensions $6$ and $3$, and $S \subset
W_{V}$ a Seifert surface of $V$ with respect to a section $s \co V \to \HV$ (namely $\bd{S} =
s(V)$) such that
$e(F(s)) = 0$ in $H^{2}(V;\Q)$.  
By \fullref{cor:seifert-e}, $2t_{S}$ is an $e$-class of $(W,V)$.
In this subsection, we prove the geometric formula for $\s(W,V,2t_{S})$.

The vanishing $e(F(s)) = 0$ implies $e(\nu_{S})|_{\bd{S}} = 0$ by \eqref{eq:Y-bdry}, and the integral
\begin{equation}
  \label{eq:int_S}
  \int_{S}^{}e(\nu_{S})^{2} \in \Q
\end{equation}
is well--defined.
(More precisely, \eqref{eq:int_S} means the integral $\int_{S}^{}a^{2}$, where $a \in
H^{2}(\pair{S};\Q)$ is an element such that the homomorphism $H^{2}(\pair{S};\Q) \to
H^{2}(S;\Q)$ maps $a$ to $e(\nu_{S})$.)

The following is the geometric formula for $\s(W,V,2t_{S})$.
\begin{theorem}
  \label{thm:geometric-formula}
  Let $(W,V)$ be a pair of closed manifolds of dimensions $6$ and $3$, and $S$ a Seifert surface
  of $V$ with respect to a section $s \co V \to \HV$ such that $e(F(s)) = 0$ over
  $\Q$.  Then, we have
  \begin{equation}
	\label{eq:geometric-formula}
	\s(W,V,2 t_{S}) = \Sign S - \int_{S}^{}e(\nu_{S})^{2}.
  \end{equation}
\end{theorem}
\begin{remark}
  The formula holds only for embedded Seifert surfaces, and not for immersed Seifert surfaces.
\end{remark}
\begin{proof}
  Set $Z = [-1,1] \x W$ which is a $7$--manifold, and 
  let us consider the submanifolds $X, Y \subset Z$ of dimensions $4$, $5$ defined by
  \begin{align*}
	X &= \left( [0,1] \x V \right) \cup_{\left\{ 0 \right\} \x V} \left( \left\{ 0 \right\} \x
	S \right)\\
	Y &= [0,1] \x S \subset Z.
  \end{align*}
  Here, we shall assume that $X$ is a smooth proper $4$--submanifold such that $\bd{X} =
  \left\{ 1 \right\} \x V$, after ``non--smooth part'' $\left\{ 0 \right\} \x V$ is rounded in a standard
  fashion, and that $Y$ has the smooth boundary $\bd{Y} = S \cup (-X)$ and the
  corner $\left\{ 1 \right\} \x V$.  
  Thus, $Y$ is a Seifert surface of $X$ in $Z$ with respect to the section $\Ts \co X \to \HX$
  such that $\Ts(X) = \HX \cap Y$. 
  By \fullref{prop:seifert},  $2t_{Y} \in H^{2}(Z \setminus X;\Q)$ is a quasi $e$-class of
  $(Z,X)$, and $e(F(\Ts))/2$ is the self--linking form of the quasi $e$-manifold
  $(Z,X,2 t_{Y})$.  By the construction, the boundary of $(Z,X,2t_{Y})$  is
  \begin{equation*}
	\bd{(Z,X,2t_{Y})} \cong (W,V,2t_{S}) \amalg (-(W, \es, 0)).
  \end{equation*}

  The oriented cobordism group $\O_{6}$ vanishes (see \eqref{tbl:O-pt}),  and so $(W,\es,0)$ bounds
  an $e$-manifold of the form $(Z',\es,0)$.  That implies $\s(W,\es,0) = 0$ by
  \fullref{axiom2}.
  Since $X \cong S$, we have $\Sign X = \Sign S$.
  By \fullref{thm:main1}\,\eqref{item:sl} and by the definition of $\L$, we have
  \begin{equation*}
	\s(W,V,2 t_{S}) = 
	\Sign X - 4 \int_{X}^{}\frac{e(F(\Ts))^{2}}{4} =
	\Sign S - \int_{X}^{}e(F(\Ts))^{2}.
  \end{equation*}
  By the Stokes' theorem, 
  \begin{equation*}
    \int_{X}^{}e(F(\Ts))^{2} =
	\int_{S}^{}e(\nu_{S})^{2} 
	- \int_{\bd{Y}}^{}e(\nu_{Y})^{2}
	=
	\int_{S}^{}e(\nu_{S})^{2},
  \end{equation*}
  where note that $\nu_{Y}|_{S} = \nu_{S}$ and $\nu_{Y}|_{X}=F(\Ts)$.
  Hence, the formula \eqref{eq:geometric-formula} holds.
\end{proof}

\subsection{\texorpdfstring{Takase's invariant}{Takase's invariant}}
\label{sec:takase}
Let $M$ be an integral homology $3$--sphere, and $f \co M \to S^{6}$ a smooth embedding.
We write $M_{f} = f(M)$ as before.
Since $(S^{6}, M_{f})$ is simple by \fullref{prop:simple}, the
invariant $\s(S^{6}, M_{f}) \in \Q$ is well--defined by \fullref{thm:Q}.
It is not difficult to show that there is a Seifert surface $S$ of $M_{f}$ with respect to some section
$s \co M_{f} \to \HM_{f}$
(cf.~\cite[Proposition\,2.5]{takase-hom3}).
Since $e(F(s)) = 0 \in H^{2}(M_{f};\Q) = 0$, the cohomology class $2t_{S}$ is the unique $e$-class of
$(S^{6},M_{f})$ by \fullref{cor:seifert-e}, and so $\s(S^{6},M_{f}) = \s(S^{6}, M_{f},
2t_{S})$.
By \fullref{thm:geometric-formula}, the geometric formula
\begin{equation*}
  \s(S^{6}, M_{f}) = \Sign S - \int_{S}^{}e(\nu_{S})^{2}
\end{equation*}
holds.
The right--hand side is nothing but ($-8$ times) the definition of Takase's invariant $\O(f)$
\cite[Proposition\,4.1]{takase-hom3},
and thus, we obtain the following immediate corollary.
\begin{corollary}
  \label{cor:takase}
  For a smooth embedding $f \co M \to S^{6}$ of an integral homology $3$--sphere $M$, we have
  \begin{equation*}
	\s(S^{6},M_{f}) = -8 \O(f).
  \end{equation*}
\end{corollary}
Since $\O(f) = H(f)$ when $M = S^{3}$ \cite[Corollary\,6.5]{takase-geom},
and so again we obtain \fullref{thm:Hf} as direct consequence of \fullref{cor:takase}.

%%%%%%%%%%%%%%%%%%%%%%%%%%%%%%%%%%%%%%%%%%%%%%%%%%%%%%%%%%%%%%%%%%%%%%%%%%%%%%%%%%%%%%%%%%%%%%%%%
\section{Milnor's triple linking number}
\label{sec:milnor}
%%%%%%%%%%%%%%%%%%%%%%%%%%%%%%%%%%%%%%%%%%%%%%%%%%%%%%%%%%%%%%%%%%%%%%%%%%%%%%%%%%%%%%%%%%%%%%%%%
In this section, we prove \fullref{thm:milnor}.
We begin by reviewing the definition of the triple linking number $\mu(L) \in \Z$ of oriented
algebraically split $3$--component links $L = K_{1} \cup K_{2} \cup K_{3}$ in $\R^{3}$ by using
Seifert surfaces.  
%-------------------------------------------------------------------------------------------
\subsection{Review of the triple linking number}
\label{sec:review-m}
%-------------------------------------------------------------------------------------------
The letters $i$ and $j$ will denote elements in $\left\{ 1,2,3 \right\}$.
Since the linking number $\lk(K_{i},K_{j})$ vanishes ($i\neq j$), 
$K_{i}$ has a Seifert surface $\S_{i}' \subset \R^{3}$, $\bd{\S_{i}'} = K_{i}$, such that
$\S'_{i} \cap K_{j} = \es$ ($i \neq j$).  
The triple linking number $\mu(L)$ is defined to be
the algebraic intersection number
$\mu(L) = \#\left( \S'_{1} \cap \S'_{2} \cap \S'_{3} \right)$.
In other words, regarding the intersection
$C_{i,j}' = \S_{i}' \cap \S_{j}'$ ($i<j$) as an oriented
$1$--dimensional closed submanifold of $\S_{i}'$, we can write
\[
\mu(L) = \#(C'_{1,2} \cap C'_{1,3}).
\]

For the proof of \fullref{thm:milnor},
we introduce slightly different (but essentially the same) definition of $\mu(L)$ as the
following.  Let $L_{0} = K_{1,0} \cup
K_{2,0} \cup K_{3,0}$ be a $3$--component unlink in $\R^{3}$ split from $L$.
Then the link $L_{i} = K_{i} \cup (-K_{i,0})$ has a connected Seifert surface $\S_{i} \subset
\R^{3}$, $\bd{\S_{i}} = L_{i}$, such that 
$\S_{i} \cap L_{j} = \es$ ($i\neq j$), and
$\mu(L)$ is defined to be
\begin{equation}\label{eq:def-mu2}
  \mu(L) = \#(C_{1,2} \cap C_{1,3}),
\end{equation}
where $C_{i,j} = \S_{i} \cap \S_{j} \subset \S_{i}$ ($i < j$)
which is an oriented  $1$--dimensional closed submanifold of $\S_{i}$.

From now on, we regard $L$ and $L_{0}$ as links in $S^{3} = \R^{3} \cup \left\{ \infty \right\}$.

%-------------------------------------------------------------------------------------------
\subsection{\texorpdfstring{Seifert surface of $M_{L}$}{Seifert surface of ML}}
\label{sec:seifert}
%-------------------------------------------------------------------------------------------
We construct an immersed Seifert surface $S$ of $M_{L}$ as follows.
First of all, let us recall the definition  of $3$--submanifold
$M_{L}$ of $T^{3} \x S^{3}$:
\begin{align*}
  T^{3}_{i} &=
  \left\{ (t_{1},t_{2},t_{3},x) \in T^{3} \x S^{3}\ | \ f_{i}(t_{i}) = x \right\}, &
  \calL &= T^{3}_{1} \cup T^{3}_{2} \cup T^{3}_{3},\\
  T^{3}_{i,0} &=
  \left\{ (t_{1},t_{2},t_{3},x) \in T^{3} \x S^{3} \ | \ f_{i,0}(t_{i}) = x \right\},&
  \calL_{0} &= T^{3}_{1,0} \cup T^{3}_{2,0} \cup T^{3}_{3,0},\\
  M_{L} &= \calL \cup (-\calL_{0}).
\end{align*}
Here, $f_{i} \co S^{1} \to S^{3}$ and $f_{i,0} \co S^{1} \to S^{3}$ are smooth embeddings
representing $K_{i}$ and $K_{i,0}$ respectively, and
$T^{3}$ is the $3$--torus with coordinates $(t_{1}, t_{2},t_{3})$ such that
$t_{i} \in S^{1} = \R/\Z$.
For convenience, we also write
\begin{equation*}
  \calL_{i} = T^{3}_{i} \cup (-T^{3}_{i,0})
\end{equation*}
so that $M_{L} = \calL_{1} \cup \calL_{2} \cup \calL_{3}$.

Since $\S_{i}$ is connected, there exists a map
\begin{equation*}
  p_{i} \co \S_{i} \to S^{1},
\end{equation*}
such that $p_{i}\, f_{i} = p_{i}\, f_{i,0} = \text{identity} \co S^{1} \to S^{1}$.
Let us consider the following smooth embeddings:
\begin{equation*}
  \begin{split}
	F_{1} \co \S_{1} \x S^{1} \x S^{1} \to T^{3} \x S^{3},\\
	F_{2} \co S^{1} \x \S_{2} \x S^{1} \to T^{3} \x S^{3},\\
	F_{3} \co S^{1} \x S^{1} \x \S_{3} \to T^{3} \x S^{3},
  \end{split}
  \qquad
  \begin{split}
	(x,t_{2},t_{3}) \mapsto
	(p_{1}(x), t_{2},t_{3}, x)\\
	(t_{1},x,t_{3}) \mapsto
	(t_{1}, p_{2}(x),t_{3}, x)\\
	(t_{1},t_{2},x) \mapsto
	(t_{1}, t_{2},p_{3}(x), x)\\
  \end{split}
\end{equation*}
The image $S_{i} \subset T^{3} \x S^{3}$ of $F_{i}$ is a Seifert surface of $\calL_{i}$, 
such that
\begin{equation}
  \label{eq:seifert}
  S_{i} \cap \calL_{j} = \es
  \quad\text{($i \neq j$)},
\end{equation}
and the union 
\begin{equation}
  \label{eq:S123}
  S = S_{1} \cup S_{3} \cup S_{3}
\end{equation}
is an immersed Seifert surface of $M_{L}$
($S_{i}$ may intersects the other components $S_{j}$ ($i \neq j$)).

The intersection
\begin{equation*}
  \S_{i,j} = S_{i} \cap S_{j} \subset S_{i}\quad\text{($i \neq j$)},
\end{equation*}
is a $\S_{i,j}$ as a $2$--dimensional closed submanifold of $S_{i}$, and
the intersection number $\#(\S_{1,2} \cap \S_{1,3}) \in \Z$ is defined.
The following lemma will be used to prove \fullref{thm:milnor} in \fullref{sec:proof6}
\begin{lemma}
  \label{lem:Smu}
  $\#(\S_{1,2} \cap \S_{1,3}) = \mu(L)$.
\end{lemma}
\begin{proof}
  Let $h \co S_{1} \to \S_{1} \x S^{1} \x S^{1}$ be the diffeomorphism defined by
  \begin{equation*}
	h(p_{1}(x), t_{2}, t_{3}, x) =
	(x, t_{2} - p_{2}(x), t_{3} - p_{3}(x))
  \end{equation*}
  for $(p_{1}(x), t_{1}, t_{3}, x) \in S_{1}$, then we have
  \begin{align*}
	h(\S_{1,2}) &= C_{1,2} \x \left\{ 0 \right\} \x S^{1},\\
	h(\S_{1,3}) &= C_{1,3} \x S^{1} \x \left\{ 0 \right\},
  \end{align*}
  and therefore, 
  \begin{align*}
	  \#(\S_{1,2} \cap \S_{1,3})
	  &=
	  \#( (C_{1,2} \x \left\{ 0 \right\} \x S^{1})
	  \cap (C_{1,3} \x S^{1} \x \left\{ 0 \right\}))\\
	  &=
	  \#(C_{1,2} \cap C_{1,3})\\
	  &=
	  \mu(L)
  \end{align*}
  by the definition \eqref{eq:def-mu2}.
\end{proof}

%-------------------------------------------------------------------------------------------
\subsection{\texorpdfstring{$(T^{3} \x S^{3}, M_{L})$ is simple}{(T3 x S3, ML) is simple}}
\label{sec:simple}
%-------------------------------------------------------------------------------------------
By using the immersed Seifert surface $S$ of $M_{L}$ constructed in the previous subsection, we prove
that $(T^{3} \x S^{3}, M_{L})$ is simple.
\begin{lemma}\label{lem:normal}
  The normal bundle $\nu_{S_{i}}$ of $S_{i} \subset T^{3} \x S^{3}$ is trivial.
\end{lemma}
\begin{proof}
  The vector field $\partial/\partial t_{i}$ on $T^{3} \x S^{3}$ is transverse to
  the submanifold $S_{i}$,
  and this gives a non--vanishing section of $\nu_{S_{i}}$.  
  Since the rank of $\nu_{S_{i}}$ is $2$, $\nu_{S_{i}}$ is trivial.
\end{proof}
Let
$t_{S} \in H^{2}( (T^{3} \x S^{3}) \setminus M_{L}; \Q)$
be the fundamental cohomology class of $S$, and we define
\begin{equation*}
  e_{L} = 2t_{S}.
\end{equation*}
\begin{proposition}\label{prop:link-e}
  The manifold pair $(T^{3} \x S^{3}, M_{L})$ has an $e$-class $e_{L}$, and is simple.
\end{proposition}
\begin{proof}
  Since the normal bundle $\nu_{S_{i}}$ is trivial by \fullref{lem:normal}, 
  the Euler class $e(F(s)) \in H^{2}(M_{L};\Q)$ vanishes, 
  where $s \co M_{L} \to \HM_{L}$ is the section such that $s(M_{L}) = \HM_{L} \cap S$.
  Thus, $e_{L}$ is an $e$-class of $(T^{3} \x S^{3}, M_{L})$ by \fullref{cor:seifert-e}.
  Since the restriction $H^{2}(T^{3}\x S^{3};\Q) \to H^{2}(M_{L};\Q)$ is injective, $(T^{3} \x S^{3},
  M_{L})$ is simple.
\end{proof}
By \fullref{thm:Q} and \fullref{prop:link-e}, we can define the invariant
\begin{equation*}
  \s(T^{3}\x S^{3}, M_{L}) = \s(T^{3} \x S^{3}, M_{L}, e_{L}) \in \Q
\end{equation*}
of $(T^{3} \x S^{3}, M_{L})$.
As we explained in \fullref{rem:link-homotopy}, this is a link homotopy invariant of $L$.

%-------------------------------------------------------------------------------------------
\subsection{\texorpdfstring{Proof of \fullref{thm:milnor}}{Proof of Theorem \ref{thm:milnor}}}
\label{sec:proof6}
%-------------------------------------------------------------------------------------------
In this subsection, we prove \fullref{thm:milnor} by using the formula in
\fullref{thm:main1}\,\eqref{item:sl}.

\begin{figure}[htbp]
  \begin{center}
	%\input{push.pstex_t}
	%%%%%%%%%%%%%%%%%%%%%%%%%%%%%%%%%%%%%%%%%%%%%%%%%%%%%%%%%%%%%%%%%%%%%%%%%%%%%

	\begin{picture}(0,0)%
	  \includegraphics{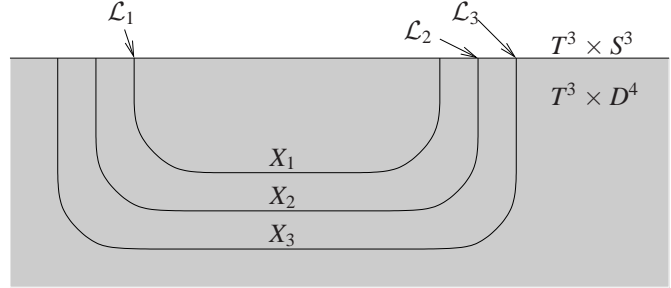}%
	\end{picture}%
	\setlength{\unitlength}{3158sp}%
	\begingroup\makeatletter\ifx\SetFigFontNFSS\undefined%
	\gdef\SetFigFontNFSS#1#2#3#4#5{%
	\reset@font\fontsize{#1}{#2pt}%
	\fontfamily{#3}\fontseries{#4}\fontshape{#5}%
	\selectfont}%
	\fi\endgroup%
	\begin{picture}(5199,2295)(-2186,-973)
	  \put(2052,476){\makebox(0,0)[lb]{\smash{{\SetFigFontNFSS{10}{12.0}{\rmdefault}{\mddefault}{\updefault}{$T^3
	  \times D^4$}%
	  }}}}
	  \put(2055,876){\makebox(0,0)[lb]{\smash{{\SetFigFontNFSS{10}{12.0}{\rmdefault}{\mddefault}{\updefault}{$T^3
	  \times S^3$}%
	  }}}}
	  \put(-149,-307){\makebox(0,0)[lb]{\smash{{\SetFigFontNFSS{10}{12.0}{\rmdefault}{\mddefault}{\updefault}{$X_2$}%
	  }}}}
	  \put(-149,
	  -4){\makebox(0,0)[lb]{\smash{{\SetFigFontNFSS{10}{12.0}{\rmdefault}{\mddefault}{\updefault}{$X_1$}%
	  }}}}
	  \put(-152,-604){\makebox(0,0)[lb]{\smash{{\SetFigFontNFSS{10}{12.0}{\rmdefault}{\mddefault}{\updefault}{$X_3$}%
	  }}}}
	  \put(901,989){\makebox(0,0)[lb]{\smash{{\SetFigFontNFSS{10}{12.0}{\rmdefault}{\mddefault}{\updefault}{$\mathcal{L}_{2}$}%
	  }}}}
	  \put(1276,1139){\makebox(0,0)[lb]{\smash{{\SetFigFontNFSS{10}{12.0}{\rmdefault}{\mddefault}{\updefault}{$\mathcal{L}_{3}$}%
	  }}}}
	  \put(-1424,1139){\makebox(0,0)[lb]{\smash{{\SetFigFontNFSS{10}{12.0}{\rmdefault}{\mddefault}{\updefault}{$\mathcal{L}_{1}$}%
	  }}}}
	\end{picture}%
	%%%%%%%%%%%%%%%%%%%%%%%%%%%%%%%%%%%%%%%%%%%%%%%%%%%%%%%%%%%%%%%%%%%%%%%%%%%%%
	\caption{Submanifold $X = X_{1} \cup X_{2} \cup X_{3}$}
	\label{fig:push}
  \end{center}
\end{figure}

We begin by constructing a proper $4$--submanifold $X \subset T^{3} \x D^{6}$
such that $\bd{X} = M_{L}$ and $X \cong S_{1} \amalg S_{2} \amalg S_{3}$.
Pushing $\Int S_{i}$ into the inside of $T^{3} \x
D^{4}$ ($\bd{S_{i}}$ is fixed on the boundary $T^{3} \x S^{3}$),
we obtain a proper $4$--submanifold
$X_{i} \subset T^{3} \x D^{4}$ such that $X_{i} \cong S_{i}$ and $\bd{X}_{i} = \calL_{i}$,
and we assume that the depth of $X_{i}$ is shallower than $X_{i+1}$ so that
$X_{i} \cap X_{j} = \es$ if $i \neq j$, see
\fullref{fig:push}.
We then define 
\[
	X = X_{1} \cup X_{2} \cup X_{3}.
\]
The natural isotopy of sinking $S_{i}$ down onto $X_{i}$ gives a $5$--submanifold
\[
	Y_{i} \subset T^{3} \x D^{4}
\]
with the boundary $\bd{Y_{i}} = S_{i} \cup (-X_{i})$ and the corner
$\calL_{i}$, and it is a Seifert surface of $X_{i}$ in $T^{3} \x D^{4}$ with respect to the section
\begin{equation*}
  s_{i} \co X_{i} \to \HX_{i}
\end{equation*}
such that  $s_{i}(X_{i}) = Y_{i} \cap \HX_{i}$.
Let $t_{Y_{i}} \in H^{2}( (T^{3} \x D^{4}) \setminus X_{i};\Q)$ be the fundamental cohomology class of
$Y_{i}$.
We define a cohomology class $\Te \in H^{2}( (T^{3} \x D^{4}) \setminus X;\Q)$ by
\begin{align*}
  \Te &= 2 (t_{Y_{1}} + t_{Y_{2}} + t_{Y_{3}}).
\end{align*}
Note that 
\begin{equation}
  \label{eq:Te-bdry}
  \Te|_{(T^{3} \x S^{3}) \setminus M_{L}} = e_{L}
\end{equation}
by definition.
We will see soon that $\Te$ is a quasi $e$-class of $(T^{3} \x D^{4}, X)$
(\fullref{prop:quasi-Y}).
\begin{remark}
The union
$
Y = Y_{1} \cup Y_{2} \cup Y_{3}
$
may not be an immersed Seifert surface of $X \subset T^{3} \x D^{4}$ in our sense,
because $X_{i} \cap Y_{j}$ may not be empty if $i < j$.
Thus, we cannot apply \fullref{prop:seifert} to $Y$ to prove that $\Te$ is a quasi $e$-class.
\end{remark}
  Let us consider the intersection
\begin{equation*}
  \S_{i,j}' = X_{i} \cap Y_{j} \subset X_{i}\quad \text{($i < j$)}
\end{equation*}
which is an oriented $2$--submanifold of $X_{i}$ such that
\begin{gather}
  \label{eq:S2} \S_{i,j}' \subset \Int X_{i},\\
  \label{eq:S3} \S_{1,2}' \cap \S_{2,3}' = \S_{1,3}' \cap \S_{2,3}' = \es.
\end{gather}

Let $\nu_{Y_{i}}$ and $\nu_{\S_{i,j}'}$ be the normal bundles of $Y_{i} \subset T^{3} \x D^{4}$ and
$\S_{i,j}' \subset X_{i}$ respectively.
\begin{lemma}\label{lem:normal-g}
  The normal bundles $\nu_{Y_{i}}$ and $\nu_{\S_{i,j}'}$ are trivial.
\end{lemma}
\begin{proof}
  The triviality of $\nu_{Y_{i}}$ follows from the definition of $Y_{i}$ and \fullref{lem:normal}.
  Since $\nu_{\S_{i,j}'}$ is isomorphic to $\nu_{Y_{j}}|_{\S_{i,j}'}$ (where we regard $\S_{i,j}' \subset
  Y_{j}$), this is trivial too.
\end{proof}
\begin{lemma}
  \label{lem:tYi}
  The cohomology class 
  $2 t_{Y_{i}}$ is an $e$--class of $(T^{3} \x D^{4}, X_{i})$.
\end{lemma}
\begin{proof}
  Since $\nu_{Y_{i}}$ is trivial by \fullref{lem:normal-g}, the cohomology class $e(F(s_{i})) \in H^{2}(X_{i};\Q)$
  vanishes.  Thus, $2t_{Y_{i}}$ is an $e$-class by \fullref{cor:seifert-e}.
\end{proof}
Let $\g_{i,j} \in H^{2}(X_{i};\Q)$ be the Poincar\'e dual
of $\S_{i,j}$, and we write
\begin{equation*}
  \g = \g_{1,2} + \g_{1,3} + \g_{2,3} \in H^{2}(X;\Q).
\end{equation*}
By the definitions of $X_{i}$ and $Y_{j}$, we have $\HX_{i} \cap Y_{j} =
\rho_{X_{i}}^{-1}(\S_{i,j}')$ ($i<j$), and this implies
\begin{equation}
  t_{Y_{j}}|_{\HX_{i}} = \rho_{X_{i}}^{*}\g_{i,j}\quad \text{($i < j$)}
  \label{eq:XY}
\end{equation}
in $H^{2}(\HX_{i};\Q)$ by the Poincar\'e duality.

We obtain the following  proposition.
\begin{proposition}\label{prop:quasi-Y}
  The triple $(T^{3} \x D^4, X, \Te)$ is a quasi $e$-manifold with the boundary $(T^{3} \x S^{3}, M_{L},
  e_{L})$ and with the self--linking form $\g$.
\end{proposition}
\begin{proof}
  By \fullref{lem:tYi} and \eqref{eq:XY},  we have
%  \[
%  \Te|_{\HX} - \sum_{i} e(F_{X_{i}})  =
%  2 (t_{Y_{2}}|_{\HX_{1}} 
%  + t_{Y_{3}}|_{\HX_{1}} 
%  + t_{Y_{3}}|_{\HX_{2}} )
%  \in H^{2}(\HX;\Q)
%  \] 
%  is the Poincar\'e dual of
%  \begin{equation*}
%	2 \sum_{i < j}
%	[\rho_{X_{i}}^{-1}(\S_{i,j})]
%	\in H_{4}(\pair{\HX};\Q).
%  \end{equation*}
%  Hence, we have
  \begin{align*}
	\Te|_{\HX} &= \sum_{i=1}^{3} e(F_{X_{i}}) + 
	2 (\rho_{X_{1}}^*\g_{1,2} + \rho_{X_{1}}^*\g_{1,3} + \rho_{X_{2}}^*\g_{2,3})\\
	&=
	e(F_{X}) +  2 \rho_{X}^*\g
  \end{align*}
  in $H^{2}(\HX;\Q)$.
  It follows from \eqref{eq:S2} that $\g|_{\bd{X}} = 0$.  Thus, $\Te$ is a quasi $e$-class of
  $(T^{3} \x D^{4}, X)$ with the self--linking form $ \g$.
  By \eqref{eq:Te-bdry}, we have $\bd{(T^{3} \x D^{4}, X, \Te)} = (T^{3} \x S^{3}, M_{L}, e_{L})$.
\end{proof}
This is the proof of \fullref{thm:milnor}.
\begin{proof}[Proof of \fullref{thm:milnor}]
  By \fullref{thm:main1}\,\eqref{item:sl} and \fullref{prop:quasi-Y},  we have
  given as follows.
  \begin{equation*}
	\s(T^{3} \x S^{3}, M_{L}) =
	\Sign X - 4 \int_{X}^{}\g^{2}.
  \end{equation*}
  Since $X_{i} \cong S_{i} \cong \pm \S \x S^{1} \x S^{1}$,
  we have $\Sign X = 0$.
  We have $\g_{i,j}^{2} = 0$ ($i \neq j$)  by \fullref{lem:normal-g},
  and $\g_{1,j} \g_{2,3} = 0$ ($j=2,3$) by \eqref{eq:S3}.  Thus,
  \begin{equation*}
	\int_{X}^{} \g^{2}=
	  2 \int_{X}^{} \g_{1,2} \g_{1,3} = 2 \#(\S_{1,2}' \cap \S_{1,3}').
  \end{equation*}

  Since the intersection $C_{1,2,3} = Y_{1} \cap Y_{2} \cap Y_{3}$ is a $1$--dimensional oriented
  cobordism from $\S_{1,2} \cap \S_{1,3}$ to $\S_{1,2}' \cap \S_{1,3}'$, namely
  \begin{equation*}
	\bd{C_{1,2,3}} = (\S_{1,2} \cap \S_{1,3}) \amalg (- (\S_{1,2}' \cap \S_{1,3}')),
  \end{equation*}
  we have
  \[
  	\#(\S_{1,2}' \cap \S_{1,3}') = \#(\S_{1,2} \cap \S_{1,3}).
  \]
  The right--hand side equals $\mu(L)$ by \fullref{lem:Smu}, and hence, we have
  \begin{equation*}
	\s(T^{3} \x S^{3}, M_{L}) = -8 \mu(M).
	\proved
  \end{equation*}
\end{proof}

%%%%%%%%%%%% References %%%%%%%%%%%%%
%%
%<Author name> is written as Initial of Given Name, and Family Name.
%<Title> is written in roman letters.
%<Journal name> should be abbreviated according to
% the MR Serials Abbreviations List of Mathematical Reviews:
% (Abbreviations of Names of Serials; http://www.ams.org/mr-database)
%For <Pages>, use en-dash "--" between page numbers.
%%

%%%%%%%%%%%%%%%%%%%%   End of main body of article
%
%                             References
%
%   BiBTeX users uncomment the following line:
%

%\bibliographystyle{gtart}

						%\bigskip
						%%%%%%%%%%%% Authors' addresses %%%%%%%%%%%%%
						%\address{Tetsuhiro Moriyama\\
						%Institut Fourier (UMR 5582), B.P. 74\\
						%38402 Saint-Martin-d'H\`eres Cedex,\\
						%France\\
						%{$e$-mail: tetsuhir@ms.u-tokyo.ac.jp,\ \ moriyama@ujf-fourier.fr}}
						%
						\end{document}